\documentclass[12pt]{amsart}
\usepackage{amssymb}

\let\cal\mathcal
\let\frak\mathfrak
\let\Bbb\mathbb

\usepackage{eucal}
\let\euc\mathcal

\def\>{\relax\ifmmode\mskip.666667\thinmuskip\relax\else\kern.111111em\fi}
\def\<{\relax\ifmmode\mskip-.333333\thinmuskip\relax\else\kern-.0555556em\fi}
\def\?{\relax\ifmmode\mskip-.666667\thinmuskip\relax\else\kern-..111111em\fi}
\def\vsk#1>{\vskip#1\baselineskip}
\def\vv#1>{\vadjust{\vsk#1>}\ignorespaces}
\def\vvn#1>{\vadjust{\nobreak\vsk#1>\nobreak}\ignorespaces}
 \let\alb\allowbreak
\let\dsty\displaystyle \let\tsty\textstyle
\let\ssty\scriptstyle \let\sssty\scriptscriptstyle
\def\fratop{\genfrac{}{}{0pt}1}
\def\satop#1#2{\fratop{\ssty#1}{\ssty#2}}
\def\tsum{\mathop{\tsty\sum}\limits}

\def\plait#1{\par\hangindent2\parindent\indent\kern\parindent
 \llap{#1\enspace}\ignorespaces}

\let\Smallskip\smallskip
\def\smallskip{\par\Smallskip}
\let\Medskip\medskip
\def\medskip{\par\Medskip}
\let\Bigskip\bigskip
\def\bigskip{\par\Bigskip}

\let\Maketitle\maketitle
\def\maketitle{\Maketitle\thispagestyle{empty}\let\maketitle\empty}

\newtheorem{thm}{Theorem}[section]
\newtheorem{cor}[thm]{Corollary}
\newtheorem{lem}[thm]{Lemma}
\newtheorem{prop}[thm]{Proposition}
\newtheorem{conj}[thm]{Conjecture}

\numberwithin{equation}{section}

\theoremstyle{definition}
\newtheorem*{rem}{Remark}
\newtheorem*{example}{Example}

\def\beq{\begin{equation}}
\def\eeq{\end{equation}}
\def\be{\begin{equation*}}
\def\ee{\end{equation*}}

\def\bean{\begin{eqnarray}}
\def\eean{\end{eqnarray}}
\def\bea{\begin{eqnarray*}}
\def\eea{\end{eqnarray*}}
\def\Ref#1{{\rm(\ref{#1})}}

\let\bt\beta
\let\gm\gamma  
\let\dl\delta \let\Dl\Delta 
 \let\eps\varepsilon \let\epsilon\eps

\let\zt\zeta
\let\tht\theta \let\Tht\Theta
\let\thi\vartheta

\let\la\lambda \let\La\Lambda
\let\rhe\varrho
\let\si\sigma 
 
 \let\phi\varphi \let\Pho\varPhi

\let\om\omega  

\let\Tilde\widetilde
\let\Hat\widehat
\let\der\partial

\let\ge\geqslant

\let\le\leqslant

\let\ox\otimes
\let\lx\ltimes

\def\nin{\not\in}

\def\C{\Bbb C}

\def\Z{\Bbb Z}

\def\gl{\frak{gl}}

\def\Bc{\cal B}

\def\Dc{\cal D}

\def\Gc{\cal G}
\def\Hc{\cal H}

\def\Mc{\cal M}
\def\Oc{\cal O}

\def\Zc{\cal Z}

\def\Ae{\euc A}
\def\Ee{\euc E}
\def\He{\euc H}
\def\Pe{\euc P}
\def\Te{\euc T}

\def\Peh{\Hat\Pe}

\def\Phh{\hat\Pho}
\def\Qh{\Hat Q}
\def\Sh{\hat S}
\def\Th{\hat T}
\def\Zh{\Hat Z}

\def\Dt{\Tilde\Dc}
\def\Het{\<\>\Tilde{\<\He\<\>}\<}
\def\Pht{\tilde\Pho}
\def\Pet{\tilde\Pe}
\def\Qt{\Tilde Q}
\def\Xt{\Tilde X}

\def\tbigoplus{\mathop{\tsty{\bigoplus}}\limits}

\def\lsym#1{#1\alb\dots\relax#1\alb} \def\lc{\lsym,}

\let\on\operatorname

\def\diag{\on{diag}}
\def\End{\on{End}}

\def\Im{\on{Im}}

\def\sign{\on{sign}}
\def\tr{\on{tr}}
\def\Tr{\on{Tr}}
\def\Wr{\on{Wr}}

\def\gln{\gl_N}
\def\glnt{\gln[t]}
\def\Ugln{U(\gln)}

\def\Uglnt{U(\glnt)}
\def\Yn{Y(\gln)}

\let\bs\boldsymbol

\let\nc\newcommand

\nc{\pone}{\C{\mathbb P}^1}

\def\chip{{\vphantom|\chi}}
\def\w{{\vrule width 0pt depth 2.6pt w}}
\def\0{{\>0}}
\def\ss{{\<\>\raise.22ex\hbox{$\sssty S$}}}
\def\yy{{\<\>\raise.22ex\hbox{$\sssty Y\!\!\!$}}}
\def\yv{{\<\>\raise.22ex\hbox{$\sssty Y$\!}}}

\def\xn{x_1\lc x_n}
\def\zn{z_1\lc z_n}
\def\pz{(\zn)}

\def\Aes{\Ae^{\?\ss}}
\def\nh{{n,\hbar}}
\def\nhb{{n,\bar\hbar}}
\def\nla{{n,\<\>\bs\la}}
\def\kml{{k,m,\<\>\bs\la}}
\def\nhla{{n,\hbar\<\>,\<\>\bs\la}}
\def\khlba{{k,\hbar\<\>,\<\>\bs\la\<\>,\bs a}}
\def\mhlba{{m,\hbar\<\>,\<\>\bs\la\<\>,\bs a}}
\def\nhlba{{n,\hbar\<\>,\<\>\bs\la\<\>,\bs a}}
\def\nnla{{n,N\<,\<\>\bs\la}}

\def\hla{{\hbar\<\>,\<\>\bs\la}}
\def\hii{{\hbar\<\>;\<\>i,i}}
\def\hij{{\hbar\<\>;\<\>i,\<\>j}}

\def\lba{{\bs\la\>,\bs a}}
\def\lbaij{{\bs\la\>,\bs a,i,\<\>j}}
\def\klba{{k,\<\>\bs\la\<\>,\bs a}}
\def\mlba{{m,\<\>\bs\la\<\>,\bs a}}
\def\nlba{{n,\<\>\bs\la\<\>,\bs a}}
\def\hlba{{\hbar\<\>,\>\bs\la\>,\bs a}}
\def\hlbami{{\hbar\<\>,\>\bs\la\>,\bs a,m,\<\>i}}
\def\rst{{\bs r\<,\<\>\bs s,\tau}}
\def\Ll{{\bs\La,\bs\la}}
\def\Llb{{\bs\La,\<\>\bs\la\<\>,\bs b}}

\def\Wrh{\Wr_\hbar}
\def\lapn{\bs\la\>\vdash\<n}
\def\lapns{{\bs\la\,\vdash\>n}}
\def\llarr{{\ssty\longleftarrow}}

\def\sing{{\it\>sing}}

\def\ia{i,\<\>a}
\def\ii{i,\<\>i}
\def\ij{i,\<\>j}

\def\ji{j,\<\>i}
\def\mi{m,\<\>i}
\def\ik{i,\<\>k}
\def\il{i,\<\>l}
\def\jk{j,\<\>k}
\def\kj{k,\<\>j}
\def\kl{k,\<\>l}
\def\%#10{#1,\<\>0}
\def\^#1{^{[#1]}}
\def\@#1{^{(#1)}}
\def\+#1{^{\{#1\}}}

\def\KZ/{{\sl KZ\/}}
\def\qKZ/{{\sl qKZ\/}}

\let\co c
\let\Co C
\let\Hh K
\def\Sg{\frak S}
\let\rr r
\let\yh h

\begin{document}

\hrule width0pt
\vsk->

\title[Bethe subalgebras of the group algebra of the symmetric group]
{Bethe subalgebras of the group algebra\\[3pt] of the symmetric group}

\author[E.\,Mukhin, V\<.\,Tarasov, and \>A.\<\,Varchenko]
{E.\,Mukhin$\>^*$, V\<.\,Tarasov$\>^\star$, and \>A.\<\,Varchenko$\>^\diamond$}

\maketitle

\begin{center}
\vsk-.2>
{\it $\kern-.4em^{*,\star}\<$Department of Mathematical Sciences,
Indiana University\,--\>Purdue University Indianapolis\kern-.4em\\
402 North Blackford St, Indianapolis, IN 46202-3216, USA\/}

\medskip
{\it $^\star\<$St.\,Petersburg Branch of Steklov Mathematical Institute\\
Fontanka 27, St.\,Petersburg, 191023, Russia\/}

\medskip
{\it $^\diamond\<$Department of Mathematics, University of North Carolina
at Chapel Hill\\ Chapel Hill, NC 27599-3250, USA\/}
\end{center}

{\let\thefootnote\relax
\footnotetext{\vsk-.8>\noindent
$^*$\,Supported in part by NSF grant DMS-0900984\\
$^\star$\,Supported in part by NSF grant DMS-0901616\\
$^\diamond$\,Supported in part by NSF grant DMS-0555327}}

\medskip
\begin{abstract}

We introduce families \,$\Bc^\ss_n\pz$ \,and \,$\Bc^\ss_\nh\pz$ \, of maximal
commutative subalgebras, called Bethe subalgebras, of the group algebra
\;$\C[\Sg_n]$ \,of the symmetric group. Bethe subalgebras are deformations
of the Gelfand\<\>-Zetlin subalgebra of \;$\C[\Sg_n]$\>.
We describe various properties of Bethe subalgebras.
\end{abstract}

\section{Introduction}
\label{intro}

Algebras of integrals of motion, or Bethe algebras, play an important role
in the theory of quantum integrable systems. There has been progress recently
in understanding properties of Bethe algebras, see \cite{MTV4}, \cite{FFRy}.

\vsk.2>
In this paper we define families \,$\Bc^\ss_n\pz$ \,and \,$\Bc^\ss_\nh\pz$
of commutative subalgebras, called Bethe subalgebras, of the group algebra
\,$\C[\Sg_n]$ of the symmetric group, and describe their properties.
Here \>$\zn$ \>and $\hbar\ne0$ \>are complex numbers. The first family is
a degeneration of the second one as \,$\hbar\to 0$. The Bethe subalgebras
of \,$\C[\Sg_n]$ \>correspond to the Bethe algebras for the Gaudin or
{\sl XXX\/}-type quantum integrable models on tensor powers of the vector
representation of the Lie algebra \,$\gl_N$ via the Schur-Weyl duality.

\vsk.2>
The subalgebras \,$\Bc^\ss_n\pz$ \,and \,$\Bc^\ss_\nh\pz$ can be viewed as
deformations of the Gelfand\<\>-Zetlin subalgebra \;$\Gc_n\<\subset\C[\Sg_n]$\>.
More precisely, for distinct \>$\zn$, the sub\-algebra \,$\Bc^\ss_n\pz$
\,is a maximal commutative subalgebra of \,$\C[\Sg_n]$ \,of dimension
\,$\dim\>\Gc_n$\>, \,and \,$\Bc^\ss_n\pz$ \,tends to \;$\Gc_n$ \,as
\vvn.3>
\beq
\label{zlim}
\frac{z_{a-1}-z_a}{z_{a+1}-z_a}\,\to\,0\,,\qquad a=2\lc n-1\,,
\vv.3>
\eeq
see Theorem~\ref{manyprop} and Proposition~\ref{GZ}. Similarly,
if \,$z_a-z_b\ne\hbar$ \>for all \;$1\le b<a\le n$, then \,$\Bc^\ss_\nh\pz$
is a maximal commutative subalgebra of \,$\C[\Sg_n]$ \,of dimension
\,$\dim\>\Gc_n$\>, \,and \,$\Bc^\ss_\nh\pz$ \,tends to \;$\Gc_n$ in
the limit~\Ref{zlim} provided in addition \,$\hbar/(z_1-z_2)\to 0$\>,
\,see Theorem~\ref{manyproph} and Proposition~\ref{GZh}.

\begin{example}
Let \>$n=3$. The center \>$\Zc_3$ \>of \,$\C[\Sg_3]$ \,is spanned over \,$\C$
\>by the identity and the elements \,$\si_{1,2}+\si_{1,3}+\si_{2,3}$\>,
\,$\si_{1,2}\>\si_{2,3}+\si_{2,3}\>\si_{1,2}$\>, where \,$\si_{a,b}$ \>is
the transposition of \,$a$ \>and \,$b$. The Gelfand\<\>-Zetlin subalgebra of
\,$\C[\Sg_3]$ \,is spanned by \>$\Zc_3$ \>and \,$\si_{1,2}$\>.
The Bethe subalgebra \,$\Bc^\ss_3(z_1,z_2,z_3)$ \,is spanned by \>$\Zc_3$
\>and the element \,$z_1\<\>\si_{2,3}+z_2\>\si_{1,3}+z_3\>\si_{1,2}$\>.
The Bethe subalgebra \,$\Bc^\ss_{3,\hbar}(z_1,z_2,z_3)$ \,is spanned
by \>$\Zc_3$ \>and the element $z_1\<\>\si_{2,3}+z_2\>\si_{1,3}+z_3\>\si_{1,2}
-\hbar\>\si_{1,2}\>\si_{2,3}$\>. Note that all maximal commutative subalgebras
of \,$\C[\Sg_3]$ \,are of this form.
\end{example}

For distinct \>$\zn$, \>consider the \KZ/ {\it elements}
\,$H_1\^n\<\lc H_n\^n\<\in\C[\Sg_n]$\>:
\vvn-.1>
\beq
\label{Ha}
H_a\^n\,=\,\sum_{\satop{b=1}{b\ne a}}^n\,\frac{\si_{a,b}}{z_a-z_b}\,,
\qquad a=1\lc n\,.
\vv-.1>
\eeq
These elements pairwise commute.
They are the right-hand sides of Knizhnik-Zamolodchi\-kov (\KZ/\>) type
equations for functions with values in $\C[\Sg_n]$, see for example~\cite{C},
\cite{FV}. If \,$\Sg_n$ acts on \,$(\C^N)^{\ox n}$ by permuting the tensor
factors, the images of \,$H_1\^n\<\lc H_n\^n$ become the right-hand sides
of the celebrated Knizhnik-Zamolodchikov equations, see~\cite{KZ}, and
the Hamiltonians of the Gaudin model associated with $\gln$, see~\cite{G}.

\vsk.2>
By Theorem~\ref{generateS}, see also Theorem~\ref{generateSt}, the elements
\,$H_1\^n\<\lc H_n\^n$ generate the Bethe subalgebra \,$\Bc^\ss_n\pz$. \,This
statement corresponds to the well-known result that the Gelfand\<\>-Zetlin
sub\-algebra of \,$\C[\Sg_n]$ \,is generated by the Young\<\>-Jucys\<\>-Murphy
elements $J_a=\sum_{b=1}^{a-1}\si_{a,b}$ \,for \,$a=2\lc n$. \,Indeed,
for every $a=2\lc n$, \,the element \>$(z_a\<-z_1)H_a\^n$ \>tends to
the Young\<\>-Jucys\<\>-Murphy element $J_a$ in the limit~\Ref{zlim}.

\vsk.2>
The counterpart of Theorem~\ref{generateS} for the Bethe subalgebra
\,$\Bc^\ss_\nh\pz$ is given by Theorem~\ref{generateSh}. It implies that
for distinct \>$\zn$ such that \,$z_a-z_b\ne\hbar$ \>for all \;$a,b=1\lc n$,
the subalgebra \,$\Bc^\ss_\nh\pz$ is generated by the \>\qKZ/ {\it elements\/}
\,$\Hh_1\^n\<\lc\Hh_n\^n$:
\vvn.2>
\begin{align}
\label{qKZ}
\Hh_a\^n\<\>={}\,
(z_a-z_{a-1}+\hbar\>\si_{a-1,\<\>a})\,\dots{} & (z_a-z_1+\hbar\>\si_{1,\<\>a})
\\[3pt]
{}\times\,(z_a-z_n+\hbar\>\si_{a,n})\dots{} &
(z_a-z_{a+1}+\hbar\>\si_{a,\<\>a+1})\,,
\notag
\\[-15pt]
\notag
\end{align}
cf.~\Ref{Ka}.
The \>\qKZ/ elements and their images in \,$\End\bigl((\C^N)^{\ox n}\bigr)$
under the action of \,$\Sg_n$ by permuting the tensor factors
are limits of the right-hand sides of difference analogues
of the \KZ/\<\>-type equations \cite{C}, \cite{FR}.

\vsk.3>
The spectra of \,$\Bc^\ss_n\pz$ \,and \,$\Bc^\ss_\nh\pz$ as commutative
algebras admit a natural description in terms of scheme-theoretic fibers
of appropriate Wronski maps, see Theorems~\ref{manyprop} and~\ref{manyproph}.
This kind of description is a refinement of the nested Bethe ansatz method
developed in the theory of quantum integrable models~\cite{KR}.

\vsk.2>
A special case of Bethe subalgebras \,$\Bc^\ss_\nh\pz$ is given by
the homogeneous Bethe subalgebra \,$\Aes_n=\Bc^\ss_\nh(z_1\lc z_1)$.
Actually, \,$\Aes_n$ \>does not depend on \,$\hbar$ \,and \,$z_1$.
The subalgebra \,$\Aes_n$ \>is a maximal commutative subalgebra of
\;$\C[\Sg_n]$ \,generated by the elements
\vvn.3>
\be
G_k\^n\,=\sum_{1\le\>i_1<\<\cdots<\>i_k\le n\!\!}
\si_{i_1,\<\>i_2}\>\si_{i_2,\<\>i_3}\dots\<\>\si_{i_{k-1},\<\>i_k}\,,
\qquad k=2\lc n\,,
\vv.1>
\ee
see Theorem~\ref{generateG}. The permutation
\,$\si_{i_1,\<\>i_2}\>\si_{i_2,\<\>i_3}\dots\<\>\si_{i_{k-1},\<\>i_k}$
is an increasing \,$k$-cycle \,$(i_1\,i_2\>\dots\>i_k)$\>.

\vsk.2>
Another nice set of generators of \,$\Aes_n$ is given by
\,$\gm_n\<=\si_{1,2}\,\si_{2,3}\>\dots\<\>\si_{n-1,n}\<=G_n\^n$ and
the {\it local charges\/} \;$I_1\^n\<\lc I_{n-2}\^n$\>, see Theorem~\ref{Inkt}
and Corollary~\ref{generateI}. The elements \;$I_k\^n$ and \;$G_k\^n$ are
related via the equality of generating series:
\,$\log\>\bigl(1+\sum_{k=1}^{n-1}\,(G_n\^n)^{-1}\>G_{n-k}\^n\>u^k\>\bigr)=
\sum_{k=1}^\infty I_k\^n\>u^k$,
see~\Ref{gmn} and~\Ref{Ik}. It is known that \;$I_1\^n\<\lc I_{n-2}\^n$ \,can
be written as sums of local densities independent of \,$n$\>, see~\Ref{Ikn},
the proof going back to~\cite{L}. In more detail, consider the chain
\;$\C[\Sg_1]\subset\C[\Sg_2]\subset\dots\subset\C[\Sg_n]\subset\cdots{}$\>,
where \;$\C[\Sg_n]$ \,is generated by the elements \,$\si_{a,a+1}$ \>for
$a=1\lc n-1$\>. Then \,for every \,$k$ \>there is an element
\;$\tht_k\?\in\Sg_{k+1}$ \,independent of \,$n$ \,such that
\vvn-.7>
\beq
\label{Ikn}
I_k\^n=\sum_{m=0}^{n-1}\,\gm_n^{\<\>m}\,\tht_k\>\gm_n^{-m}\,,
\qquad k=1\lc n-2\,,
\eeq
cf.~\Ref{Ink}; for instance, \;${\tht_1=\si_{1,2}}$\>,
\,\,${\tht_2=(\si_{2,3}\,\si_{1,2}-\si_{1,2}\,\si_{2,3}-1)/2}$\>.
\,Notice that the image of \,$I_1\^n$ in $\End\bigl((\C^2)^{\ox n}\bigr)$
under the action of \,$\Sg_n$ \,is essentially the Hamiltonian of
the celebrated {\sl XXX\/} Heisenberg model, whose eigenvectors and
eigenvalues were first studied in~\cite{B}.

\vsk.2>
The algebra \,$\Aes_n$ \>is semisimple, and its action on every irreducible
representation of \,$\Sg_n$ \,has simple spectrum, see Theorem~\ref{lessprop}.

\vsk.2>
The group algebra of the symmetric group enters various important families
of algebras. An interesting question is if there are analogues of Bethe
subalgebras for other members of those families. There are indications
that such analogues may exist for the Weyl groups of root systems other than
of type $A$, though there are several gaps to be closed there. The finite Hecke
algebras of type $A$ and Birman-Wenzl-Murakami algebras probably have versions
of Bethe subalgebras due to their relation to centralizer constructions and
analogues of the Schur-Weyl duality. The corresponding integrable models
should be the models with reflecting boundary conditions, see for
example~\cite{I}. It is plausible that the Bethe subalgebras can be defined
for group algebras of affine Weyl groups.

\vsk.2>
The algebra \,$\Bc^\ss_n\pz$ is closely related to the center
of the rational Cherednik algebra of type $A_n$ at the critical level,
see~\cite{MTV8}. We expect that the Bethe subalgebra
\,$\Bc^\ss_\nh\pz$ has a similar relation to the center of
the trigonometric Cherednik algebra. An interesting open question is
to describe an analogue of the Bethe subalgebra related to the center
of the double affine Hecke algebra.

\vsk.2>
The plan of the paper is as follows.
We introduce the Bethe subalgebras \,$\Bc^\ss_n\pz$ \,in
Section~\ref{CSn}. In Section~\ref{Gaudin}, we review  of the Gaudin model. More sophisticated
properties of \,$\Bc^\ss_n\pz$ are described in Section~\ref{more}.
In Section~\ref{A3}, we review the definition and properties of the Bethe
sublagebra of the Yangian \,$\Yn$. We define the Bethe subalgebras
\,$\Bc^\ss_\nh\pz$ in Section~\ref{xxx} and study their properties in
Section~\ref{morexxx}. In Section~\ref{homo}, we consider the homogeneous Bethe
subalgebra \,$\Aes_n$. Additional technical details are given in Appendix.

\vsk.2>
The authors thank M.\,Nazarov and A.\>Vershik for useful discussions and
referees for valuable remarks.

\vsk->\vsk0>

\section{Bethe subalgebras \,$\Bc^\ss_n\pz$ \,of \,$\C[\Sg_n]$}
\label{CSn}

Let \,$\Sg_m$ \,be the symmetric group on \,$m$ \,symbols.
For distinct $r_1\lc r_m\<\in\{\>1\lc n\>\}$ \,we denote by
\,$\pi_{r_1\lc r_m}\^n:\C[\Sg_m]\to\C[\Sg_n]$ \,the embedding induced by
the correspondence \,$i\>\mapsto r_i$.
\vsk.2>
Let \,$A\^m=\frac1{m\<\>!}\sum_{\si\in\Sg_m}(-1)^\si\si$ \,be the
antisymmetrizer in $\C[\Sg_m]$\>; in particular, \,$A\^1=1$. Given complex
numbers \,$\zn$\>, consider the polynomials \,$\Pho_1\^n(u)\lc\Pho_n\^n(u)$
\,in one var\-i\-able with coefficients in $\C[\Sg_n]$\,:
\beq
\label{Phin}
\Pho_i\^n(u)\,=\!\sum_{1\le r_1<\<\dots<r_i\le n}\!\! i\<\>!\;
\pi_{r_1\lc r_i}\^n(A\^i)\!\!\prod_{\satop{a=1}{a\nin\{r_1\lc r_i\}}}^n
\!\!\!(u-z_a)\,=\,\sum_{j=0}^{n-i}\,\Pho_{\ij}\^n\,u^{n-i-j}\,.
\vv-.1>
\eeq
For instance, \;$\Pho_1\^n(u)=\sum_{r=1}^n\prod_{a\ne r}(u-z_a)$ \;and
\vvn.3>
\beq
\label{Phoi0}
\Pho_{\%i0}\^n\,=\!\sum_{1\le r_1<\<\dots<r_i\le n}\!\!
i\<\>!\;\pi_{r_1\lc r_i}\^n(A\^i)\,,\qquad i=1\lc n\,.
\eeq

\begin{prop}[\cite{OV}]
\label{center}
The elements \;$\Pho_{\%10}\^n\>\lc\Pho_{\%n0}\^n$
\;generate the center of \,$\C[\Sg_n]$\>.
\end{prop}
\noindent
Independently, this statement follows from Proposition~\ref{Zx}.

\vsk.5>
Denote by \,$\Bc^\ss_n\pz$ \,the subalgebra of $\C[\Sg_n]$ generated
by all \,$\Pho_{\ij}\^n$, \,for \,$i=1\lc n$, \,$j=0\lc n-i$.
The subalgebra $\Bc^\ss_n\pz$ depends on $\zn$ as parameters.
\vvn.3>
Clearly,
\beq
\label{Bzz}
\Bc^\ss_n\pz\,=\,\si\,\Bc^\ss_n(z_{\si(1)}\lc z_{\si(n)})\,\si^{-1}
\vv.3>
\eeq
for any \,$\si\in\Sg_n$. We call the subalgebras \,$\Bc^\ss_n\pz$
\;{\it Bethe subalgebras\/} of \,$\C[\Sg_n]$ of \,{\it Gaudin type\/}.

\vsk.2>
Set
\vvn-.6>
\beq
\label{Phiuv}
\Pho\^n(u,v)\,=\,v^n\,\prod_{a=1}^n\,(u-z_a)\>
+\>\sum_{i=1}^n\,(-1)^i\,\Pho_i\^n(u)\,v^{n-i}\,,
\vv.2>
\eeq
where \,$v$ \,is an indeterminate.
\begin{lem}
\label{Phiuvil}
We have
\vvn-.2>
\beq
\label{Phiuvi}
\Pho\^n(u,v)\,=\,(-1)^n\sum_{\si\in\Sg_n}\si\,\sign(\si)
\prod_{\satop{b=1}{b=\si(b)\!\!}}^n\bigl(1-v\>(u-z_b)\bigr)\,.
\eeq
\end{lem}
\begin{proof}
For \,$\si\in\Sg_n$, let \,$I_\si\?\subset\{1\lc n\}$ \>be the set of its fixed
points. Then formula \Ref{Phin} yields
\be
\Pho\^n(u,v)\,=\,(-1)^n\sum_{\si\in\Sg_n}\si\,\sign(\si)
\sum_{J\subset I_\si}\,(-v)^{n-|J|}\,\prod_{b\<\>\in J}\,(u-z_b)\,,
\vv-.2>
\ee
which proves the claim.
\end{proof}

\begin{lem}
\label{Bsz}
We have
\;$\Bc^\ss_n(sz_1\lc sz_n)=\Bc^\ss_n\pz$ \,for any \,$s\ne 0$,
\;and \;$\Bc^\ss_n(z_1+s\lc z_n+s)=\Bc^\ss_n\pz$ \,for any \,$s$.
\end{lem}
\begin{proof}
Formula~\Ref{Phin} yields
\;$\Pho_i\^n(su;sz_1\lc sz_n)=s^{n-i}\>\Pho_i\^n(u;\zn)$.
\,Hence,
\;$\Pho_{\ij}\^n(sz_1\lc sz_n)=s^j\>\Pho_{\ij}\^n\pz$,
\,which proves the first claim. Similarly, the second claim follows from
the equality \;$\Pho_i\^n(u+s;z_1+s\lc z_n+s)=\Pho_i\^n(u;\zn)$.
\end{proof}

\begin{prop}
\label{Bcomm}
The subalgebra \,$\Bc^\ss_n\pz$ is commutative.
\end{prop}
\begin{proof}
Consider the algebra
\,$\Hc\>=\>\C[\xn]\ox\bigl(\C[y_1\lc y_n]\lx\C[\Sg_n]\bigr)$
\,and the polynomial in $u,v$
\vvn-.3>
\be
\Phh\^n(u,v)\,=\,(-1)^n
\sum_{\si\in\Sg_n}\>\si\,\sign(\si)
\prod_{b=\si(b)}\bigl(1-(u-x_b)\>(v-y_b)\bigr)\,=\,
\sum_{i,j=0}^n\,\Phh_{\ij}\^n\,u^{n-j}\>v^{n-i}
\ee
with coefficients in $\Hc$. It is shown in~\cite{MTV8} \,(see Theorem~2.5,
Lemma~3.1 and Section~3.2 therein) that the elements \,$\Phh_{\ij}\^n$
\>commute with each other.

Let \,$\zt:\Hc\to\C[\Sg_n]$ \,be the homomorphism defined by the assignment
\vvn.3>
\be
x_i\mapsto z_i\,,\quad y_i\mapsto 0\,,\quad \si\mapsto\si\,,
\qquad i=1\lc n\,,\quad\si\in\Sg_n\,.
\vv.3>
\ee
Then \,$\zt\bigl(\Phh_{\ij}\^n\bigr)=0$ \,for $j<i$ \,and
\,$\zt\bigl(\Phh_{\ij}\^n\bigr)=\Pho_{\ij\<-i}\^n$ \,for $j\ge i\ge 1$\>,
\,which yields the claim.
\end{proof}

Alternatively, Proposition~\ref{Bcomm} follows from Theorem~\ref{BnN}
and Corollary~\ref{BB} below.

\vsk.5>
Let $\Zc_m$ be the center of $\C[\Sg_m]$. The subalgebra
$\Gc_n\subset\C[\Sg_n]$ generated by the images $\pi_{1\lc m}\^n(\Zc_m)$,
\,$m=1\lc n$, is called the Gelfand\<\>-Zetlin subalgebra of $\C[\Sg_n]$.
The subalgebra $\Gc_n$ is a maximal commutative subalgebra of $\C[\Sg_n]$,
see~\cite{OV}.

\begin{prop}
\label{GZ}
The Bethe subalgebra \,$\Bc^\ss_n\pz$ \,tends to the Gelfand\<\>-Zetlin
subal\-gebra \,$\Gc_n$ \,as \,$(z_{a-1}-z_a)/(z_a-z_{a+1})\to0$ \;for all
\,$a=2\lc n-1$\>.
\end{prop}
\begin{proof}
Without loss of generality we can assume that $z_1=0$, see Lemma~\ref{Bsz},
so we have \,$z_a/z_{a+1}\to0$ \;for all \,$a=2\lc n-1$\>. In this limit
the element \,$(-1)^j\>\Pho_{\ij}\^n\,z_{n-j+1}^{-1}\dots z_n^{-1}$ tends to
$\pi_{1\lc n-j}\^n\bigl(\Pho_{\%i-j0}\^{n-j\<\>}\bigr)$. Therefore, the limit
of \,$\Bc^\ss_n\pz$ contains \,$\Gc_n$, see Proposition~\ref{center}.
Since \,$\Gc_n$ is a maximal commutative subalgebra of \,$\C[\Sg_n]$,
and \,$\Bc^\ss_n\pz$ \,is commutative for any $\zn$,
the limit of \,$\Bc^\ss_n\pz$ coincides with \,$\Gc_n$.
\end{proof}

Let $\bt\mapsto\bt^{\>\dag}$ be the linear antiinvolution on \,$\C[\Sg_n]$
\,such that \,$\si^\dag\<=\si^{-1}$ for any \,$\si\in\Sg_n$.
Let $\bt\mapsto\bt^{\>*}$ be the semilinear antiinvolution on \,$\C[\Sg_n]$
\,such that \,$\si^*\<=\si^{-1}$ for any \,$\si\in\Sg_n$.

\begin{prop}
\label{B*}
We have
\,$\bigl(\Pho\^n_{\ij}\pz\bigr)^\dag\<=\Pho\^n_{\ij}\pz$
\,and
\;$\bigl(\Pho\^n_{\ij}\pz\bigr)^{\<*}\<=
\Pho\^n_{\ij}(\bar z_1\lc\bar z_n)$ \,for all \,$i,j$. \,Here
\,$\bar z_1\lc\bar z_n$ \>are the complex conjugates of \,$\zn$\>.
\qed
\end{prop}

Further properties of the subalgebras \,$\Bc^\ss_n\pz$
are given in Section~\ref{more}.

\section{Bethe algebra of the Gaudin model}
\label{Gaudin}

Let \,$V=\C^N$. We identify elements of \,$\End(V)$ \,with \,$N\,{\times}\,N$
\,complex matrices. We also consider $V$ as the natural vector representation
of the group \,$GL_N$.

Let \,$E_{\ij}\<\in\End(V)$ \,be the matrix with only one nonzero entry equal
to $1$ at the intersection of the $i$-th row and $j$-th column. Consider
first-order differential operators in $u$ with \,$\End(V^{\ox n})$-valued
coefficients:
\be
X_{\ij}\,=\,\dl_{\ij}\,\der_u\>-\>
\sum_{a=1}^n\,\frac{1^{\ox(a-1)}\ox E_{\ij}\ox 1^{\ox(n-a)}}{u-z_a}\;,
\qquad i,j=1\lc N\,,
\vv.2>
\ee
where \,$\dl_{\ij}$ is the Kronecker symbol,
and the $n$-th order differential operator in $u$
\vvn.2>
\be
\Dc\,=\,\prod_{a=1}^n\,(u-z_a)\,\sum_{\si\in\Sg_N}\sign(\si)\,
X_{\si(1),1}\>X_{\si(2),2}\dots X_{\si(N),N}\,.
\vv.2>
\ee
By Theorem~1 in~\cite{MTV6}, \;$\Dc$ \,is a polynomial differential operator,
\vvn.1>
\be
\Dc\,=\,\sum_{i=0}^n\,\sum_{j=0}^{n-i}\,
(-1)^i\,\Co_{\ij}\^n\,u^{\<\>n-i-j}\>\der_u^{\>N\?-i}\,.
\ee
Denote by \,$\Bc_{n,N}\pz$ \,the subalgebra of \,$\End(V^{\ox n})$
\,generated by all \,$\Co_{\ij}\^n$ \;for \,$i=1\lc n$, \,$j=0\lc n-i$.
The algebra \,$\Bc_{n,N}\pz$ depends on $\zn$ as parameters.

\begin{thm}[\cite{T}]
\label{BnN}
The algebra \,$\Bc_{n,N}\pz$ \>is commutative and commutes with
the action of \;$GL_N\?$ on \>$V^{\ox n}$.
\end{thm}

We call the algebra \,$\Bc_{n,N}\pz$ the {\it Bethe algebra
for the Gaudin model\/} with parameters $\zn$.

\begin{rem}
The algebra \,$\Bc_{n,N}\pz$ is the image of the Bethe subalgebra
of $\Uglnt$ in the tensor product $\ox_{a=1}^nV(z_a)$ of evaluation
$\glnt$-modules. We give more details in~Section~\ref{A1}.
\end{rem}

Let the symmetric group \,$\Sg_n$ act naturally on $V^{\ox n}$ by permuting
the tensor factors. Denote by \,$\varpi_n:\C[\Sg_n]\to\End(V^{\ox n})$
\,the corresponding homomorphism. By Theorem~\ref{BnN}, the algebra
\,$\Bc_{n,N}\pz$ lies in the image of \,$\varpi_n$ \>due to
the Schur-Weyl duality.

\begin{thm}
\label{SW}
We have \,$\varpi_n\bigl(\Bc^\ss_n\pz\bigr)=\Bc_{n,N}\pz$.
More precisely, \,$\varpi_n(\Pho_{\ij}\^n)=\Co_{\ij}\^n$ \,for any
\,$i=1\lc n$, \,$j=0\lc n-i$.
\end{thm}
\begin{proof}
Let \,$v_1\lc v_N$ be the standard basis of \,$V$. Identify the space
\,$V^{\ox n}\<$ with the span of monomials \,$x_{i_1,1}\dots x_{i_n,n}$
in extra variables \,$x_{\ij}$ \,by sending each vector
\,$v_{i_1}\!\lsym\ox v_{i_n}$ to the monomial \,$x_{i_1,1}\dots x_{i_n,n}$.
Under this identification the matrix of the operator
\,$x_{\ia}\>\der_{x_{j,a}}$ equals \,$1^{\ox(a-1)}\ox E_{\ij}\ox 1^{\ox(n-a)}$
\>and \;$\varpi_n\bigl(k\<\>!\;\pi_{r_1\lc r_k}\^n(A\^k)\bigr)$ \>is the matrix
of the operator
\vvn.3>
\be
\sum_{1\le a_1<\<\dots<a_k\le n}\!\!
\det(x_{r_s,\<\>a_t})_{s,\<\>t=1}^k\,\det(\der_{x_{r_s,a_t}})_{s,\<\>t=1}^k\,.
\vv.2>
\ee
Thus the claim follows from Theorem~1 in~\cite{MTV6}.
\end{proof}

\begin{cor}
\label{BB}
The algebra \,$\Bc_{n,N}\pz$ \,for \,$N\ge n$ \,is isomorphic to
\,$\Bc^\ss_n\pz$.
\end{cor}
\begin{proof}
The homomorphism \,$\varpi_n:\C[\Sg_n]\to\End(V^{\ox n})$ \,is injective for
\,$N\ge n$.
\end{proof}

A partition \,$\bs\la=(\la_1,\la_2,\ldots{})$ \,with at most \>$m$ \>parts is
a sequence of integers
such that \,$\la_1\ge\la_2\ge\dots\ge 0$ \,and \,$\la_{m+1}\<=0$\>.
If \,$\sum_{i=1}^\infty\la_i=n$\>, \,we write \,$\lapn$
\>and say that \,$\bs\la$ \,is a partition of \>$n$\>.

\vsk.2>
For \,$\lapn$, \,let \,$M_{\bs\la}$ \>be the irreducible \,$\Sg_n$-module
corresponding to \,$\bs\la$\>, \,and \,$\chip_{\bs\la}\<\in\Zc_n$
\,the respective central idempotent. By definition, \,$\chip_{\bs\la}$ acts
as the identity on $M_{\bs\la}$ and acts
as zero on any $M_{\bs\mu}$ for $\bs\mu\ne\bs\la$.
By Proposition~\ref{center}, the algebra \,$\Bc^\ss_n\pz$ \,contains \,$\Zc_n$,
so \,$\chip_{\bs\la}\?\in\Bc^\ss_n\pz$.
Set \,$\Bc^\ss_\nla\pz=\chip_{\bs\la}\Bc^\ss_n\pz$.
The algebra \,$\Bc^\ss_\nla\pz$ \,is isomorphic to the image of
\vvn.1>
\,$\Bc^\ss_n\pz$ in \,$\End(M_{\bs\la})$ by the canonical projection.
\,Clearly, \,$\Bc^\ss_n\pz=\tbigoplus_\lapns\Bc^\ss_\nla\pz$.

\vsk.2>
By the Schur-Weyl duality, we have the decomposition
\beq
\label{VLM}
V^{\ox n}\,=\,\bigoplus_{\satop\lapns {\la_{N+1}=\,0}}
L_{\bs\la}\ox M_{\bs\la}
\eeq
with respect to the \,$GL_N\<\times\Sg_n$ \,action. Here \,$L_{\bs\la}$ \,is
the irreducible representation of \,$GL_N$ with highest weight $\bs\la$.
By Theorem~\ref{BnN}, the action of \,$\Bc_{n,N}\pz$
on \,$V^{\ox n}$ descends to the action on each \,$M_{\bs\la}$.
Denote by \,$\Bc_\nnla\pz$ the image of \,$\Bc_{n,N}\pz$
\,in \,$\End(M_{\bs\la})$.

\begin{cor}
\label{Bnl}
Let \,$\lapn$ be such that \,$\la_i=0$ \>for \,$i>N$.
Then the algebra \,$\Bc_\nnla\pz$ is isomorphic to
\,$\Bc^\ss_\nla\pz$.
\end{cor}
\begin{proof}
The claim follows from Theorem~\ref{SW}.
\end{proof}

\begin{prop}
\label{Zx}
We have
\be
\sum_{i=0}^n\,(-1)^i\,\Pho_{\%i0}\^n\prod_{j=i+1}^n(t+j)\,=\,
\sum_\lapns\,\chip_{\bs\la}\,\prod_{j=1}^n\,(t-\la_j+j)\,,
\ee
where \,$t$ is an indeterminate, \>$\Pho_{0,\<\>0}\^n=1$,
\>other \>$\Pho_{\%i0}\^n$ \>are given by~\Ref{Phoi0},
\>and \,$\bs\la=(\la_1,\la_2,\ldots{})$.
\end{prop}
\begin{proof}
\vsk.2>
Let \,$e_{\ij}=\sum_{a=1}^n 1^{\ox(a-1)}\ox E_{\ij}\ox 1^{\ox(n-a)}\in
\End(V^{\ox n})$\>. \,Set
\be
(V^{\ox n})_{\bs\la}^\sing\,=\,\{\,v\in V^{\ox n}\ |
\ e_{\ii}\>v\>=\>\la_i\>v\,,
\ \ e_{j,k}\>v=0\,,\quad i=1\lc N\,,\ \ 1\le j<k\le N\,\}\,.
\ee
The subspace $(V^{\ox n})_{\bs\la}^\sing$ is nonzero if and only if
\,$\la_i=0$ \>for \>$i>N$. \,By the Schur-Weyl duality, a nonzero subspace
$(V^{\ox n})_{\bs\la}^\sing$ is an $\Sg_n$-submodule of \,$V^{\ox n}$
isomorphic to $M_{\bs\la}$. Now the proposition follows from Theorem~\ref{SW},
Corollary~\ref{Bnl}, and the results of~\cite{MTV4},
see formulae~(2.11) and~(2.3) therein.
\end{proof}

\section{Further properties of the Bethe subalgebras \,$\Bc^\ss_n\pz$}
\label{more}

The algebras \>$\Bc_\nnla$ have been studied in~\cite{MTV4}, see more details
in~Section~\ref{A1}. Theorems~\ref{BFrob} and~\ref{manyprop} translate
the results of~\cite{MTV4} into properties of the algebras $\Bc^\ss_n$
using Corollaries~\ref{BB} and~\ref{Bnl}.

\begin{thm}
\label{BFrob}
We have
\begin{enumerate}
\item[i)]
For any $\zn$, \,the algebra \,$\Bc^\ss_n\pz$ is a Frobenius
algebra.
\item[ii)]
For real \,$\zn$, the algebra \,$\Bc^\ss_n\pz$ is a direct sum
of the one-dimensional algebras isomorphic to \,$\C$. This assertion holds
for generic complex \,$\zn$ as well.
\qed
\end{enumerate}
\end{thm}

We refer a reader to~\cite{Wik} for the definition and basic properties of
Frobenius algebras.

\vsk.5>
Let \,$\Sg_n$ act on \,$\C[y_1\lc y_n]$ \,by permuting the variables.
Denote by \,$\deg\>p$ \,the homogeneous degree of \,$p\in\C[y_1\lc y_n]$\,:
\vvn.1>
\,$\deg\>y_i=1$ \,for all \,$i=1\lc n$\>. \,We extend the degree
to \,$M_{\bs\la}\ox\alb\C[y_1\lc y_n]$ \,trivially on the first factor.
\vvn.16>
Then the \,$\Sg_n$-module $M_{\bs\la}\ox\alb\C[y_1\lc y_n]$ is graded.
For any \,$w\in(M_{\bs\la}\ox\alb\C[y_1\lc y_n])^{\Sg_n}$
\vvn.16>
we have \,$\deg\>w\ge\sum_{i=1}^n(i-1)\>\la_i$, and the component of
\,$(M_{\bs\la}\ox\alb\C[y_1\lc y_n])^{\Sg_n}$ of degree
\vvn.16>
$\sum_{i=1}^n(i-1)\>\la_i$ \>is one-dimensional, see~\cite{K}.
Let \;$w_{\bs\la}$ \>be a nonzero element of
\,$(M_{\bs\la}\ox\alb\C[y_1\lc y_n])^{\Sg_n}$
of degree $\sum_{i=1}^n(i-1)\>\la_i$.

\vsk.2>
For a positive integer \,$m$ \,and a partition \,$\bs\la$ \,with at most \>$m$
\>parts, consider indeterminates \>$f_{\ij}$ with \,$i=1\lc m$ \,and
\,$j=1\lc\la_i+m-i$, \;$j\ne\la_i-\la_s-i+s$ \,for \,$s=i+1\lc m$.
Given in addition a collection of complex numbers \,$\bs a=(a_1\lc a_n)$,
define the algebra $\Oc_\mlba$ as the quotient of
\,$\C[f_{1,1}\lc f_{m,\la_m}]$
\,by relations~\Ref{WOla} described below.

\vsk.2>
Consider $f_{\ij}$ as coefficients of polynomials in one variable,
\beq
\label{fiu}
f_i(u)\,=\,u^{\la_i+m-i}+\sum_{j=1}^{\la_i+m-i} f_{\ij}\,u^{\la_i+m-i-j}\,,
\eeq
assuming that \,$f_{i,\la_i-\la_s-i+s}=0$ \,for \,$s>i$,
\vvn.1>
that is,
the coefficient of \,$u^{\la_s+m-s}$ \,in \>$f_i(u)$ \>equals zero.
The defining relations for \,$\Oc_\mlba$ \,are written as an equality of
two polynomials in \>$u$:
\beq
\label{WOla}
\Wr[f_1(u)\lc f_m(u)]\,=\prod_{1\le i<j\le m}\!(\la_j-\la_i+i-j)\;
\Bigl(u^n+\sum_{s=1}^n\,(-1)^s\,a_su^{n-s}\Bigr)\,,
\vv.3>
\eeq
where
\;$\Wr[f_1(u)\lc f_m(u)]=\det\<\>\bigl(\der_u^{\,i-1}f_j(u)\bigr)_{i,j=1\lc m}$
\;is the Wronskian.

\begin{lem}
\label{Okm}
Let \,$\bs\la$ \,be a partition with at most \>$m$ \>parts, and \,$k\ge m$.
Then the algebras \,$\Oc_\klba$ \>and \,$\Oc_\mlba$ \>are isomorphic.
\end{lem}
\begin{proof}
Let \,$f_{\ij}^{\{k\}}$ \>and \,$f_{\ij}^{\{m\}}$ \>be the indeterminates used
\vvn.1>
to define the algebras \,$\Oc_\klba$ \>and \,$\Oc_\mlba$\>, respectively.
In both cases, the subscripts \,${i\<\>,j}$ \,run through the same sets
of pairs because \,$\la_i=0$ \>for $i>m$\>. \,Denote by \,$f_i^{\{k\}}(u)$
\,and \,$f_i^{\{m\}}(u)$ the corresponding polynomials, see~\Ref{fiu}.
Notice that \,$f_i^{\{k\}}(u)=u^{k-i}$ \,for \,$i=m+1\lc k$.

\vsk.2>
The assignment
\vvn-.5>
\be
f_{\ij}^{\{m\}}\,\mapsto\,
f_{\ij}^{\{k\}}\prod_{s=m+1}^k\!\frac{\la_i+s-i-j}{\la_i+s-i}
\vv.3>
\ee
defines an isomorphism of \,$\Oc_\mlba$ \>and \,$\Oc_\klba$.
\,Indeed, the assignment means that
\beq
\label{fi}
f_i^{\{m\}}(u)\,\mapsto\,
\der_u^{\>k-m}f_i^{\{k\}}(u)\,\prod_{s=m+1}^k\,\frac1{\la_i+s-i}
\vv-.3>
\eeq
\vvn-.4>
and
\begin{gather*}
\Wr[f_1^{\{m\}}(u)\lc f_m^{\{m\}}(u)]\,\mapsto\,
C_\kml\>\Wr[f_1^{\{k\}}(u)\lc f_k^{\{k\}}(u)]\,,
\\[4pt]
C_\kml\>=\,(-1)^{(k-m)(k+m-1)/2}\>\prod_{s=0}^{k-m-1}\frac1{s\<\>!}
\,\;\prod_{i=1}^m\,\frac{(\la_i+m-i)\<\>!}{(\la_i+k-i)\<\>!}\;,
\end{gather*}
which yields the claim.
\end{proof}

\vsk.2>
Let \,$\Oc_\lba=\Oc_\nlba$. \,Set
\be
F_\lba(u,v)\,=\,e^{-uv}\>\Wr[f_1(u)\lc f_n(u),e^{uv}\>]\>
\prod_{1\le i<j\le m}\frac1{\la_j-\la_i+i-j}\;.
\ee
It is a polynomial in \,$u,v$ \,with coefficients in \,$\Oc_\lba$\>:
\vvn.1>
\be
F_\lba(u,v)\,=\,\sum_{i=0}^n\,
\sum_{j=0}^{n-i}\,(-1)^i\,F_\lbaij\,u^{\<\>n-i-j}\>v^{\<\>n-i}\,.
\vv.2>
\ee

Further on, we identify elements of \,${M_{\bs\la}\ox\alb\C[y_1\lc y_n]}$
with \,$M_{\bs\la}$-valued polynomials in \>$y_1\lc y_n$\>.

\begin{thm}
\label{manyprop}
Let \,$\zn$ be distinct. Then
\begin{enumerate}
\item[i)]
The algebra \,$\Bc^\ss_n\pz$ is a maximal commutative subalgebra of
\,$\C[\Sg_n]$.
\vsk.2>
\item[ii)]
The map \,$\Bc^\ss_n\pz\to\tbigoplus_\lapns M_{\bs\la}$\>,
\;$X\mapsto\tbigoplus_\lapns Xw_{\bs\la}\pz$\>, is an isomorphism
of the regular representation of \,$\Bc^\ss_n\pz$ on itself
\vvn.1>
and the \,$\Bc^\ss_n\pz$-module $\tbigoplus_\lapns M_{\bs\la}$\>.
In particular, \,$\dim\Bc^\ss_n\pz=\sum_\lapns\dim M_{\bs\la}$\;.
\vsk.2>
\item[iii)]
If \,$\zn$ are real, then the action of \,$\Bc^\ss_n\pz$ on
\,$\tbigoplus_\lapns M_{\bs\la}$ is diagonalizable and has simple
spectrum. This assertion holds for generic complex \,$\zn$ as well.
\vsk.2>
\item[iv)]
The assignment \,$\chip_{\bs\la}\Pho_{\ij}\^n\mapsto F_\lbaij$ \,for
\,$i=1\lc n$, \,$j=0\lc n-i$, extends to an isomorphism of algebras
\,$\Bc^\ss_\nla\pz\to\Oc_\lba$.
Here \,$\bs a=(a_1\lc a_n)$ \>and
\;$u^n+\sum_{s=1}^n\,(-1)^s\,a_su^{n-s}=\prod_{i=1}^n(u-z_i)$\,.
\end{enumerate}
\end{thm}
\begin{proof}[Proof of Theorems~\ref{BFrob}, \ref{manyprop}]
Due to the decomposition $\Bc^\ss_n\pz=\tbigoplus_\lapns\Bc^\ss_\nla\pz$,
\,it suffices to verify the counterparts of the claims for the algebras
$\Bc^\ss_\nla\pz$ and the $\Bc^\ss_\nla\pz$-modules $M_{\bs\la}$. The required
statements follow from the properties of the algebras \,$\Bc_\nnla\pz$\>,
established in~\cite{MTV4}, see Sections~5.2 and 5.3 \>{\it op.\;sit\/}.
The correspondence of notation between this paper and~\cite{MTV4} is described
in Section~\ref{A1}.
\end{proof}

\begin{rem}
If $\zn$ coincide at most in pairs, the algebra \,$\Bc^\ss_n\pz$ \,is a maximal
commutative subalgebra of \,$\C[\Sg_n]$. If there is a triple of coinciding
$z$'s, then \,$\Bc^\ss_n\pz$ \,is not a maximal commutative subalgebra of
\,$\C[\Sg_n]$.
\end{rem}

Let \,$\si_{a,b}\in\Sg_n$ \>denote the transposition of \,$a$ \>and \,$b$.
\>For distinct \,$\zn$, \,set
\be
H_a\^n\,=\,\sum_{\satop{b=1}{b\ne a}}^n\,\frac{\si_{a,b}}{z_a-z_b}\,,
\qquad a=1\lc n\,,
\vv-.3>
\ee
cf.~\Ref{Ha}. It is easy to see that
\vvn-.2>
\be
\Pho_2\^n(u)\,=\,\sum_{a=1}^n\,\,\Bigl(
-\>H_a\^n\>+\>\sum_{b\ne a} \frac1{z_a-z_b}\,\Bigr)
\prod_{\satop{b=1}{b\ne a}}^n\,(u-z_b)\,.
\vv-.5>
\ee

Consider the diagonal matrix
\vvn.2>
\beq
\label{Z}
Z\,=\,\diag\pz
\vv.2>
\eeq
and the matrix
\vvn-.6>
\beq
\label{Q}
Q\;=\,\left(\,
\begin{matrix}
\yh_1 & \dfrac{1}{z_1-z_2} & \dfrac{1}{z_1-z_3} & \,\dots & \dfrac{1}{z_1-z_n}
\\[14pt]
\dfrac{1}{z_2-z_1} & \yh_2 & \dfrac{1}{z_2-z_3} & \,\dots & \dfrac{1}{z_2-z_n}
\\[9pt]
{}\dots & {}\dots & {}\dots & \,\dots & {}\dots
\\[6pt]
\dfrac{1}{z_n-z_1} & \dfrac{1}{z_n-z_2}& \dfrac{1}{z_n-z_3}& \,\dots & \yh_n
\end{matrix}\,\right)
\vv.4>
\eeq
depending on new variables $\yh_1\lc\yh_n$. Set
\vvn.4>
\beq
\label{Pe}
\Pe(u,v;\zn;h_1\lc h_n)\,=\,\det\<\>\bigl((u-Z)\>(v-Q)-1\bigr)\,.
\vv.2>
\eeq

\begin{thm}
\label{generateS}
Let \,$\zn$ be distinct. Then the subalgebra \,$\Bc^\ss_n\pz$
is generated by the elements \,$H_1\^n\<\lc H_n\^n$. More precisely
\vvn.3>
\beq
\label{PhiPe}
\Pho\^n(u,v)\,=\,\Pe(u,v;\zn;H_1\^n\<\lc H_n\^n)\,,
\vv.3>
\eeq
where \,$\Pho\^n(u,v)$ \,is given by~\Ref{Phiuv}.
\end{thm}
\begin{proof}
The claim follows from Theorem~\ref{SW} and Corollary~\ref{BB},
and~\cite[Theorem~3.2\>]{MTV7}. Note that the matrix \>$Q$ \>here is transposed
compared with its counterpart in~\cite{MTV7}.
\end{proof}

Combining formulae~\Ref{Phiuvi} and~\Ref{PhiPe}, we get
\vvn.3>
\be
\Pe(u,v;\zn;H_1\^n\<\lc H_n\^n)\,=\,(-1)^n
\sum_{\si\in\Sg_n}\,\si\,\sign(\si)\prod_{b=\si(b)}\bigl(1-v\>(u-z_b)\bigr)\,.
\ee

\vsk-.1>
In the limit \,${(z_{a-1}\<-z_a)\big/(z_a\<-z_{a+1})\to0}$ \,for
\,$a=2\lc n-1$\>, \,the elements \>${(z_b\<-z_1)H_b\^n}\<$ \,tend to
the Young\<\>-Jucys\<\>-Murphy elements \;$J_b=\sum_{a=1}^{b-1}\si_{a,b}$\>.
The elements \,$J_2\lc J_n$ \>gen\-er\-ate the Gelfand\<\>-Zetlin subalgebra
$\Gc_n$, see~\cite{OV}. Theorem~\ref{generateS} is the counterpart of this fact
for the Bethe subalgebra $\Bc^\ss_n\pz$.

\vsk.2>
The subalgebra of $\C[\Sg_n]$, generated by the elements
\,$H_1\^n\<\lc H_n\^n$, \>and its relation to
the Gelfand\<\>-Zetlin subalgebra were considered in~\cite{CFR}.

\vsk.2>
For distinct $\zn$ and a partition $\bs\la=(\la_1\lc\la_n)$ of $n$,
define the algebra \,$\He_{\bs\la}\pz$ as the quotient of
\,$\C[h_1\lc h_n]$ \,by relations~\Ref{relH}, \Ref{relHl} described below.
Write
\be
\Pe(u,v;\zn;h_1\lc h_n)\,=\,
\sum_{i,j=0}^n\,\Pe_{\ij}(\zn;h_1\lc h_n)\,u^{n-j}\>v^{n-i}\,.
\ee
The defining relations for \,$\He_{\bs\la}\pz$ are
\vvn.3>
\beq
\label{relH}
\Pe_{\ij}(\zn;h_1\lc h_n)\,=\,0\,,\qquad 0\le j<i\le n\,,
\eeq
and
\vvn-.5>
\beq
\label{relHl}
\sum_{i=0}^n\,\Pe_{\ii}(\zn;h_1\lc h_n)\prod_{j=i+1}^n(t+j)\,=\,
\prod_{j=1}^n\,(t-\la_j+j)\,,
\vv.2>
\eeq
where \,$t$ is a formal variable, cf.~Proposition~\ref{Zx}.

\begin{thm}
\label{BSnH}
Let \,$\zn$ be distinct. Then the assignment
\,$h_a\mapsto\chip_{\bs\la}H_a\^n$, \,$a=1\lc n$, \,defines an isomorphism
of algebras \;$\He_{\bs\la}\pz$ \,and \,$\Bc^\ss_\nla\pz$\>.
\end{thm}

\begin{proof}
Define the algebra \,$\Ee_{\bs\la}\pz$ as the quotient of
\,$\C[r_1\lc r_n]$
\,by relations~\Ref{relE}, \Ref{relEl}, described below.

\vsk.2>
Let $g_a(v) = (v+\rr_a)\>e^{z_av}$, \,$a=1\lc n$.
\ Let \ $\dsty\Dl\,=\,\prod_{1\le a<b\le n}(z_b-z_a)$\,. \ Write
\vvn-.1>
\be
\Wr\bigl[g_1(v)\lc g_n(v), e^{uv}\bigr]\,=\,\Dl\,e^{uv+\sum_{a=1}^n z_av}\,
\sum_{i,j=0}^n\,R_{\ij}(\zn,\rr_1\lc\rr_n)\,u^{n-j}\>v^{n-i}\,.
\ee
The defining relations for \,$\Ee_{\bs\la}\pz$ \,are
\vvn.3>
\beq
\label{relE}
R_{\ij}(\zn,\rr_1\lc\rr_n)\,=\,0\,,\qquad 0\le j<i\le n\,,
\eeq
and
\vvn-.5>
\beq
\label{relEl}
\sum_{i=0}^n\,R_{\>\ii}(\zn,\rr_1\lc\rr_n)\prod_{j=i+1}^n(t+j)\,=\,
\prod_{j=1}^n\,(t-\la_j+j)\,.
\vv.2>
\eeq
Formulae~(4.1)\,--\,(4.3) in~\cite{MTV7} and formulae
\Ref{relH}\,--\,\Ref{relEl} in this paper imply that the algebras
$\Ee_{\bs\la}\pz$ and $\He_{\bs\la}\pz$ are isomorphic
by the correspondence
\vvn.3>
\be
h_a\,\mapsto\,-\>\rr_a-\sum_{b\ne a}\,\frac1{z_a-z_b}\;,\qquad a=1\lc n\,.
\ee
The algebra \,$\Ee_{\bs\la}\pz$ \,was studied
in~\cite[Section~5]{MTV2}, see details in Section~\ref{A2}.
Theorem~\ref{BSnH} follows from Theorem~6.12 in~\cite{MTV2}, Theorem~3.1
in~\cite{MTV6}, Lemma~4.3 in~\cite{MTV7}, and Corollary~\ref{Bnl}
in this paper.
\end{proof}

For a partition \,$\bs\la$ \,of \,$n$\>, set
\;$\dsty{\Pi_{\bs\la}(v)\,=\,\prod_{i=1}^n\,\prod_{j=1}^{\la_i}\,(v+i-j)}$\;.
\vvn-.8>
\;Notice that \;$\Pi_{\bs\la}(v)$ \,is a polynomial in \,$v$ \>of degree \,$n$
\>and
\vvn.2>
\beq
\label{Pila}
\frac{\Pi_{\bs\la}(v)}{\Pi_{\bs\la}(v+1)}\,=\,
\prod_{i=1}^n\,\frac{v-\la_i+i}{v+i}\;.
\eeq
\,$\Pi_{\bs\la}(v)$ \>is the product over all boxes of the Young
diagram of shape \>$\bs\la$\>, the number \>$i-j$ \>in the factor \>$v+i-j$
\,being the negative of content of the \>$(i,j)$-th box of the Young diagram.
\vsk.2>
Set \,\;$\dsty\Pi(v)\,=\,\sum_\lapns\>\chip_{\bs\la}\>\Pi_{\bs\la}(v)$\;.
\,It is known that
\vvn-.4>
\beq
\label{PiJ}
\Pi(v)\,=\,\prod_{i=1}^n\,(v-J_i)\,=\,
(-1)^n\sum_{\si\in\Sg_n}\>\si\,(-v)^{c(\si)}\,=\,
\sum_{\si\in\Sg_n}\>\si\,\sign(\si)\,v^{\<\>c(\si)}\>,
\eeq
see~\cite{J}. Here $J_1=0$, \>$J_2\<\>\lc J_n$ are
the Young\<\>-Jucys\<\>-Murphy elements, and $c(\si)$ is
the number of \,$\si$-orbits in \,$\{\>1\lc n\>\}$\>.

\vsk.2>
Consider a new generating function
\vvn.2>
\beq
\label{Phtuv}
\Pht\^n(u,v)\,=\,\Pi(v+1)\,\biggl(\;\prod_{a=1}^n\,(u-z_a)\,+\>
\sum_{i=1}^n\,(-u)^i\,\Pho_i\^n(u)\,\prod_{j=1}^i\,\frac1{v+j}\,\biggr)\>,
\eeq
where \,$\Pho_1\^n(u)\lc\Pho_n\^n(u)$ \,are given by~\Ref{Phin},
cf.~\Ref{Phiuv}.
\begin{lem}
The function \,$\Pht(u,v)$ is a polynomial in \>$v$ \>of degree \>$n$\>.
\end{lem}
\begin{proof}
Since \,$\Pi_{\bs\la}(v+1)$ is divisible by \,$(v+1)\dots(v+i)$ \,provided
\,$\la_i\ne 0$, it suffices to show that \,$\chip_{\bs\la}\>\Pho_i\^n(u)=0$
\,if \,$\la_i=0$\>. This follows from formula~\Ref{Phin}, because
\,$\pi_{r_1\lc r_i}\^n(A\^i)$ \,acts by zero in the \,$\Sg_n$-module
\,$M_{\bs\la}$ \,if \,$\la_i=0$\>, \,that is,
\,$\chip_{\bs\la}\>\pi_{r_1\lc r_i}\^n(A\^i)=0$\>.
\end{proof}

Set
\vvn-.8>
\be
\Pht\^n(u,v)\,=\,\sum_{i=0}^n\,\Pht_i\^n(v)\,u^{n-i}\>.
\vv.3>
\ee
By Proposition~\ref{Zx}, and formulae~\Ref{Pila} and~\Ref{Phtuv}, we have
\vvn.3>
\beq
\label{Pht0n}
\Pht_0\^n(v)\,=\,\Pi(v)\,,\qquad
\Pht_n\^n(v)\,=\,(-1)^n\>z_1\dots z_n\,\Pi(v+1)\,.
\eeq

\vsk.3>
Let \,$Z,Q$ \, be the matrices given by~\Ref{Z}, \Ref{Q}. Set
\vvn.2>
\be
\Pet(u,v;\zn;h_1\lc h_n)\,=\,\det\<\>\bigl((u-Z)\>(v-ZQ)-Z\<\>\bigr)\,.
\ee
\be
\Pet_0(v;\zn;h_1\lc h_n)\,=\,\det\<\>(v-ZQ)\,.
\vv.2>
\ee

\begin{thm}
\label{generateSt}
Let \,$\zn$ be distinct. Then
\vvn.2>
\be
\Pht\^n(u,v)\,=\,\Pet(u,v;\zn;H_1\^n\<\lc H_n\^n)\,.
\vv.2>
\ee
\end{thm}
\begin{proof}
The claim follows from Theorem~\ref{SW} and Corollary~\ref{BB}, and
appropriately modified \cite[Theorem~3.2\>]{MTV7}. The proof of the required
counterpart of \cite[Theorem~3.2\>]{MTV7} is similar to the proof of
the original assertion in~\cite{MTV7}, using the results from~\cite{MTV3}.
\end{proof}

\begin{cor}
\label{PetPi}
We have \,\;$\Pet_0(v;\zn;H_1\^n\<\lc H_n\^n)\,=\,\Pi(v)$\,.
\end{cor}
\begin{proof}
The claim follows from each of equalities~\Ref{Pht0n}.
\end{proof}

In the limit \,${z_a/z_{a+1}\to0}$ \,for \,$a=1\lc n-1$\>, \,the elements
\,$z_1H_1\^n\<\lc z_nH_n\^n$ \>tend to the the Young\<\>-Jucys\<\>-Murphy
elements \,$J_1\<\>\lc J_n$\>, the matrix \>$Q$ \>becomes triangular, and
Corollary~\ref{PetPi} reproduces the first part of formula~\Ref{PiJ}.
Moreover, in this limit
\vvn.3>
\be
(-1)^i\>z_{n-\<\>i\<\>+1}^{-1}\dots z_n^{-1}\,\Pht_i\^n(v)\,\to\,
(v-J_1)\dots(v-J_{n-\<\>i})\>(v+1-J_{n-\<\>i\<\>+1})\dots(v+1-J_n)
\vv.3>
\ee
for \,$i=1\lc n$\>. Recall that \,$J_1\<=0$\>.

\vsk.3>
For a rational function \,$R(v)$ \,denote by \,$\langle R(v)\rangle$
\vvn.1>
\,the fractional part of \,$R(v)$\>. That is, $R(v)-\langle R(v)\rangle$ is
a polynomial and \;$\lim_{\,v\to\infty}\langle R(v)\rangle=0$.

\vsk.2>
For distinct $\zn$ and a partition $\bs\la=(\la_1\lc\la_n)$ of $n$,
define the algebra \,$\Het_{\bs\la}\pz$ \,as the quotient of
\,$\C[h_1\lc h_n]$ \,by relations~\Ref{relHtl}, \Ref{relHt} described below.
Write
\be
\Pet(u,v;\zn;h_1\lc h_n)\,=\,
\sum_{i=0}^n\,\Pet_i(v;\zn;h_1\lc h_n)\,u^{n-i}\,.
\ee
The defining relations for \,$\Het_{\bs\la}\pz$ are
\vvn.3>
\beq
\label{relHtl}
\Pet_0(v;\zn;h_1\lc h_n)\,=\,\Pi_{\bs\la}(v)\,,
\eeq
cf.~Corollary~\ref{PetPi}, and
\vvn-.5>
\beq
\label{relHt}
\biggl\langle\,\frac{\Pet_i(v)}{\Pi_{\bs\la}(v+1)}\;
\prod_{j=1}^{n-i}\,(v+j)\biggr\rangle\,=\,0\,,\qquad i=1\lc n\,.\kern-2em
\eeq

\begin{conj}
\label{BSnHt}
Let \,$\zn$ be distinct. Then the assignment
\,$h_a\mapsto\chip_{\bs\la}H_a\^n$, \,$a=1\lc n$, \,defines an isomorphism
of algebras \;$\Het_{\bs\la}\pz$ \,and \,$\Bc^\ss_\nla\pz$\>.
\end{conj}

\section{Bethe subalgebra of the Yangian \,$\Yn$}
\label{A3}
The Yangian $\Yn$ is a unital asso\-ciative algebra with generators
\,$t_{\ij}\+s$ \,for \,$i,j=1\lc N$, \;$s\in\Z_{>0}$, \,subject to relations
\vvn.2>
\be
(u-v)\>\bigl[\<\>t_{\ij}(u)\>,t_{\kl}(v)\<\>\bigr]\>=\,
t_{\kj}(v)\>t_{\il}(u)-t_{\kj}(u)\>t_{\il}(v)\,,\qquad i,j,k,l=1\lc N\,,
\vv.3>
\ee
where \;$t_{\ij}(u)=\dl_{\ij}+\sum_{s=1}^\infty\,t_{\ij}\+s\>u^{-s}$\>.
\vvn.16>
The Yangian $\Yn$ is a Hopf algebra with the co\-product
\;$\Dl:\Yn\to\Yn\ox\Yn$ \,given by
\vvn.16>
\;$\Dl\bigl(t_{\ij}(u)\bigr)=\sum_{k=1}^N\,t_{\kj}(u)\ox t_{\ik}(u)$ \,for
\,$i,j=1\lc N$\>. The Yangian $\Yn$ contains \>$\Ugln$ \,as a Hopf subalgebra,
the embedding given by \,$e_{\ij}\mapsto t_{\ji}\+1$.
\vsk.2>
Set
\beq
\label{Tem}
\Te_m(u)\,=\,\sum_{\si\in\Sg_m}\,\sum_{1\le\>i_1<\<\dots<\>i_m\le n\!\!}
\sign(\si)\,t_{i_{\si(1)},\<\>i_1}(u-m+1)\dots\>t_{i_{\si(m)},\<\>i_m}(u)\,.
\vv.2>
\eeq
The series \,$\Te_1(u)\lc\Te_N(u)$\>, are called {\it transfer-matrices\/},
see~\cite{KS} and a recent exposition in~\cite{MTV1}. Formula~\Ref{Tem}
is obtained from~formulae~(4.13) and~(4.9) in~\cite{MTV1}.

\vsk.2>
The subalgebra \;$\Bc^\yy_N\<\subset\Yn$ \,generated by coefficients of
the series \,$\Te_1(u)\lc\Te_N(u)$ \,is called the {\it Bethe subalgebra\/}
of \,$\Yn$\>.

\begin{thm}[\cite{KS}]
\label{BY}
The subalgebra \,$\Bc^\yy_N$ is commutative and commutes with the subalgebra
$\Ugln\subset\Yn$.
\qed
\end{thm}

Let \>$E_{\ij}\<\in\End(\C^N)$ be as in Section~\ref{Gaudin}. Given a complex
number \>$x$, the assignment \,$t_{\ij}\+s\mapsto E_{\ji}\>x^{s-1}$
\,for \,$i,j=1\lc N$, \,\;$s\in\Z_{>0}$, that is,
\vvn.1>
\,$t_{\ij}(u)\mapsto \dl_{\ij}+E_{\ji}\>(u-x)^{-1}$, \>makes \,$V\?=\<\>\C^N$
into the evaluation module \,$V^\yv(x)$ \>over \,$\Yn$.

\vsk.2>
Consider the decomposition \Ref{VLM} of \,$V^{\ox n}\<$ with respect to
the \,$GL_N\<\times\Sg_n$ \,action. By Theorem~\ref{BY}, the action of
\,$\Bc^\yy_N$ on the \,$\Yn$-module \,$V^\yv(x_1)\ox\dots\ox V^\yv(x_n)$
descends to the action on each \,$M_{\bs\la}$. We denote by
\,$\Bc^\yy_\nnla(\xn)$ the corresponding image of \,$\Bc^\yy_N$
in \,$\End(M_{\bs\la})$.

\vsk.2>
Properties of the algebras \,$\Bc^\yy_\nnla(\xn)$ \,are studied in
\cite{MTV10}. They are similar to the properties of the algebras
\,$\Bc_\nnla\pz$ \,for the Gaudin model, established in \cite{MTV4}.
The algebra \,$\Bc^\yy_\nnla(\xn)$ coincides with the algebra
\,$\Bc^{\<\>\bs q=\bs1}\bigl(V(\bs x)^{\textit{sing}}_{\bs\la}\bigr)$,
where \,$\bs x=(\xn)$, in \cite{MTV10}.

\section{Bethe subalgebras \,$\Bc^\ss_\nh\pz$ \,of \,$\C[\Sg_n]$}
\label{xxx}

In this section we describe further generalization of the Bethe subalgebras
of \,$\C[\Sg_n]$, depending on an additional nonzero complex parameter $\hbar$.

\vsk.2>
Let $p$ be a complex number. Define the linear map
\,$\Tr_m\@p:\C[\Sg_{n+m}]\to\C[\Sg_n]$ \,by the rule:
\plait{a)} write an element \,$\si\in S_{n+m}$ \>as a product of cycles, remove
all symbols $j$ such that $j>n$ from the record of $\si$, and read the obtained
record as an element $\tau\in\Sg_n$\>;
\plait{b)} let \>$c(\si)$ be the number of \,$\si$-orbits
in \,$\{\>1\lc m+n\>\}$ \,and \>$c(\tau)$ the number of
\,$\tau$-orbits in \,\,$\{\>1\lc n\>\}$\,;
\plait{c)} then, \,$\Tr_m\@p(\si)\>=\>p^{\>c(\si)-\<\>c(\tau)}\>\tau$\>.

\begin{example}
Let \,$n=4$, \,$m=5$, \,$\si=(137)\>(256)\>(89)\in\Sg_9$\>.
Then \,$\tau=(13)\<\in\Sg_4$\>, \;$c(\si)=4$\>, \;$c(\tau)=3$\>,
\,and \,$\Tr_5\@p(\si)=p\>(13)\in\C[\Sg_4]$\>.
\end{example}

\vsk.2>
The definition of the map \,$\Tr_m\@p$ \,is motivated by
Proposition~\ref{trace}. Recall that \>$V=\C^N$, and $\Sg_k$ acts on
$V^{\ox k}$ by permuting the tensor factors, the corresponding homomorphism
denoted by \,$\varpi_k:\C[\Sg_k]\to\End(V^{\ox k})$. Define the linear map
\,$\tr_m:\End(V^{\ox(n+m)})\to\End(V^{\ox n})$ as the trace over the last $m$
tensor factors, that is, $\tr_m(X\ox Y)=X\>\tr(Y)$ for any $X\<\in\End(V^{\ox
n})$ and $Y\in\End(V^{\ox m})$.

\begin{prop}
\label{trace}
We have \;$\tr_m\bigl(\varpi_{n+m}(\si)\bigr)=
\varpi_n\bigl(\Tr_m\@N(\si)\bigr)$ \,for any \,$\si\in S_{n+m}$.
\qed
\end{prop}

Recall that \,$\pi_{r_1\lc r_m}\^{n+m}:\C[\Sg_m]\to\C[\Sg_{n+m}]$ \,is
\vvn.2>
the embedding induced by the correspondence \,$i\>\mapsto r_i$.
For brevity, we set \,$\thi_m=\pi_{n+1\lc n+m}\^{n+m}$.

\begin{lem}
\label{TrXY}
Let \,$r_1\lc r_k$ \,and \,$s_1\lc s_l$ \,be two collections of numbers
from \,$\{\<\>1\lc n+m\>\}$, distinct within each collection, and
\,$\{r_1\lc r_k\<\>\}\cap\{s_1\lc s_l\<\>\}\subset\{\<\>n+1\lc n+m\>\}$\>.
Then for any \,$X\<\in\Sg_k$ \,and \;$Y\?\in\Sg_l$\>, \,we have
\vvn.3>
\beq
\label{XY}
\Tr_m\@p\bigl(\pi_{r_1\lc r_k}\^{n+m}(X)\,\pi_{s_1\lc s_l}\^{n+m}(Y)\>\bigr)
\,=\,\Tr_m\@p\bigl(\pi_{s_1\lc s_l}\^{n+m}(Y)\,
\pi_{r_1\lc r_k}\^{n+m}(X)\bigr)\,.
\vv.3>
\eeq
\end{lem}
\begin{proof}
Since the homomorphism \,$\varpi_n:\C[\Sg_n]\to\End(V^{\ox n})$ \,is injective
for $N\ge n$, formula \Ref{XY} holds for \,$p\in\Z_{\ge n}$
\,by Proposition~\ref{trace} and properties of the trace of matrices.
This implies formula~\Ref{XY} \,for all \,$p$ \,because both sides depend
on \,$p$ \,polynomially,
\end{proof}

Recall that \,$A\^m\in\C[\Sg_m]$ \,is the antisymmetrizer,
\,$A\^m=\frac1{m\<\>!}\sum_{\si\in\Sg_m}(-1)^\si\si$\,.
\begin{lem}
\label{TrA}
Let \,$k\le m$, \,and \,$X\in\Sg_{n+k}$. Then
\be
\Tr_m\@p\bigl(\thi_m(A\^m)\>\pi_{1\lc n+k}\^{n+m}(\<X)\bigr)\,=\,
\Tr_k\@p\bigl(\thi_k(A\^k)\>X\bigr)\,\prod_{i=1}^{m-k}\,\frac{p-m+i}{m+1-i}\;.
\ee
\end{lem}
\begin{proof}
Straightforward by induction on \>$m$.
\end{proof}

Recall that \,$\si_{a,b}$ \>is the transposition of \,$a$ \>and \,$b$.
For \,$m\in\Z_{\ge 0}$, consider the polynomials \,$T_m\^n(u;p;\hbar\>)$
\,in \>$u$ \>with coefficients in $\C[\Sg_n]$ depending on \>$p$ \>and
\>$\hbar$ \>as parameters:
\begin{align}
\label{T0}
T_0\^n(u;p;\hbar\>)\, &{}=\,\prod_{a=1}^n\,(u-z_a)\,,
\\[2pt]
\label{Tnm}
T_m\^n(u;p;\hbar\>)\,&{}=\,\Tr_m\@p\Bigl(\thi_m(A\^m)\>
\bigl(u-z_n+\hbar\tsum_{i=1}^m\si_{n,n+i}\bigr)\dots
\bigl(u-z_1+\hbar\tsum_{i=1}^m\si_{1,n+i}\bigr)\<\Bigr)\,.
\end{align}
Set
\vvn-.4>
\beq
\label{Tmi}
T_m\^n(u;p;\hbar\>)\,=\,\sum_{i=0}^nT_{\mi}\^n(\<\>p;\hbar\>)\,u^{n-i}\,.
\vv.2>
\eeq
Denote by \,$\Bc^\ss_\nh\pz$ \,the unital subalgebra of \,$\C[\Sg_n]$
\,generated by all \;$T_{\mi}\^n(\<\>p;\hbar\>)$ \,for \,$m=1\lc n-1$,
\;$i=1\lc n$, \,with given \,$p\>,\hbar$\>.
The subalgebra $\Bc^\ss_\nh\pz$ depends on \,$\zn$ as parameters.
We call the subalgebras \,$\Bc^\ss_\nh\pz$ \;{\it Bethe subalgebras\/}
of \,$\C[\Sg_n]$ \>of \,{\sl XXX\/} \,{\it type\/}.

\begin{prop}
\label{Bp}
The subalgebra \,$\Bc^\ss_\nh\pz$ \,does not depend on \,$p$.
\vsk-.1>
\end{prop}
\begin{proof}
We construct below elements \,$S_{\mi}\^n\<\in\C[\Sg_n]$ \,that do not depend
\vvn-.16>
on \>$p$\>, see~\Ref{Snmi}, and show in Proposition~\ref{BS} that
\,$\Bc^\ss_\nh\pz$ \,is generated by \,$S_{\mi}\^n$ \,for \;$m=1\lc n-1$,
\;$i=1\lc n-m$.
\end{proof}

\begin{lem}
\label{Xlem}
We have
\beq
\label{sX}
\bigl(u-z_n+\hbar\tsum_{i=1}^m\si_{n,n+i}\bigr)\dots
\bigl(u-z_1+\hbar\tsum_{i=1}^m\si_{1,n+i}\bigr)\,=\!
\sum_{k=0}^{\min\>(m,n)}\!
\sum_{\satop{\satop{1\le r_1<\<\dots<r_k\le n}{1\le s_1<\<\dots<s_k\le m\!\!}}
{\tau\in\Sg_k\!\!}}\!\!\hbar^{\<\>k}\>X_\rst(u;\hbar\>)\,,
\vv-.7>
\eeq
where
\vvn-.5>
\begin{align}
\label{X}
X_\rst(u;\hbar\>)\,={}& \!\prod_{r_k\<<\<\>b\<\>\le n}^\llarr\!
\bigl(u-z_a+\hbar\tsum_{j=1}^k\si_{b,n+s_j}\>\bigr)\,\times{}
\\[3pt]
{}\times{}&\prod_{1\<\>\le a\le k}^\llarr
\Bigl(\si_{r_{\<a},n+s_{\tau(a)}}
\prod_{r_{\<a-1}\<<\<\>b\<\><r_{\<a}}^\llarr\!
\bigl(u-z_b+\hbar\tsum_{j=1}^{a-1}\si_{b,n+s_{\tau(j)}}\>\bigr)\<\Bigr)\,,
\kern-1.6em
\notag
\\[-16pt]
\notag
\end{align}
$\dsty\prod^\llarr$ \,denotes the ordered product:
$\dsty\prod_{c\>\le\>i\<\>\le\<\>d}^\llarr x_i\>=\,x_d\,x_{d-1}\dots x_c$\,,
\;and \;$r_0\<=0$.
\end{lem}
\begin{proof}
Expand both sides of~\Ref{sX} as a sum of monomials
\,$y_n\,y_{n-1}\dots y_1$\>, where each factor \>$y_b$ \>is either \,$(u-z_b)$
\,or \,$\si_{b,n+j}$ \>for some \>$j\in\{1\lc m\}$. Then one can verify
by inspection that every of \,$(m+1)^n$ possible monomials appear exactly once
in each side of~\Ref{sX}.
\end{proof}

Set \;$S_0\^n(u;\hbar\>)=\prod_{a=1}^n(u-z_a)$,
\beq
\label{Sk}
S_k\^n(u;\hbar\>)\,=\,k\<\>!\;\hbar^{\<\>k}\!\!
\sum_{1\le r_1<\<\dots<r_k\le n}\!\!\pi_{r_1\lc r_k}\^n(A\^k)\!\!
\prod_{\satop{1\le\<\>a\<\>\le\<\>n}{a\nin\{r_1\lc r_k\}}}^\llarr\!\!\!
\bigl(u-z_a+\hbar\tsum_{\satop{j=1}{r_{\<j}<a}}^k\si_{r_{\<j},\<\>a}\>\bigr)\,.
\eeq
for \,$k=1\lc n$, \,and \;$S_k\^n(u;\hbar\>)=0$ \;for \,$k>n$.
\,For example, since \,$A\^1=1$\>,
\vvn.3>
\begin{align}
\label{S1}
S_1\^n(u;\hbar\>)\, &{}=\,\hbar\,\sum_{i=1}^n\,
(u-z_n+\hbar\>\si_{i,n})\dots(u-z_{i+1}+\hbar\>\si_{\ii+1})\>
(u-z_{i-1})\dots(u-z_1)
\\[2pt]
& {}=\,\sum_{j=1}^n\,\sum_{1\le\>i_1<\<\cdots<\>i_j\le n\!\!}\hbar^{\>j}\>
\si_{i_1,\<\>i_j}\>\si_{i_1,\<\>i_{j-1}}\dots\<\>\si_{i_1,\<\>i_2}
\prod_{a\nin\{i_1\lc i_j\}}\!(u-z_a)\,,
\notag
\\[-14pt]
\notag
\end{align}
where the empty product of transpositions for \>${j=1}$ \>is the identity.
Notice that the permutation
\,$\si_{i_1,\<\>i_j}\>\si_{i_1,\<\>i_{j-1}}\dots\<\>\si_{i_1,\<\>i_2}\<=
\si_{i_1,\<\>i_2}\>\si_{i_2,\<\>i_3}\dots\<\>\si_{i_{j-1},\<\>i_j}$
\,is an increasing \,$j$-cycle \,$(i_1\,i_2\>\dots\>i_j)$\>.

\begin{prop}
\label{TS}
For \,$m\in\Z_{\ge 0}$\>, \>we have
\vvn.1>
\be
T_m\^n(u;p;\hbar)\,=\,\sum_{k=0}^m\,\frac1{(m-k)\<\>!}\;S_k\^n(u;\hbar)\,
\prod_{i=1}^{m-k}\,(p-m+i)\,.
\vv.3>
\ee
\end{prop}
\begin{proof}
By formula~\Ref{Sk} and Lemma~\ref{XY}, we have
\vvn.5>
\be
S_k\^n(u;\hbar)\,=\,k\<\>!\;\hbar^{\<\>k}\!\!\sum_{1\le r_1<\<\dots<r_k\le n}
\!\!\Tr_k\@p\bigl(\thi_k(A\^k)\>X_\rst(u;\hbar\>)\bigr)\,,
\vv.3>
\ee
where \,$X_\rst(u;\hbar\>)$ \,is given by~\Ref{X}, and the equality holds
for every \,$\bs s$ \,and \>$\tau$. Then the claim follows by
formulae~\Ref{Tnm}, \Ref{sX}, and Lemma~\ref{TrA}.
\end{proof}

Consider the generating series
\vvn.1>
\beq
\label{ThSh}
\Th\^n(u,x\<\>;p;\hbar)\,=\,\sum_{m=0}^\infty T_m\^n(u;p;\hbar)\,x^m\,,\qquad
\Sh\^n(u,x;\hbar)\,=\,\sum_{m=0}^\infty S_m\^n(u;\hbar)\,x^m\,.
\eeq
Remind that \,$\Sh(u,x;\hbar)$ \>is actually a polynomial in \>$x$ \>of degree
\vvn.1>
\>$n$, because \,$S_m\^n(u;\hbar)=0$ \,for \>$m>n$. Proposition~\ref{TS} is
equivalent to
\;$\Th(u,x\<\>;p;\hbar)=(1+x)^p\,\Sh\bigl(u,x/(1+x);\hbar\bigr)$\>.
Therefore, \;$\Sh(u,x;\hbar)=(1-x)^p\,\Th\bigl(u,x/(1-x);p;\hbar\bigr)$\>,
which gives
\vvn.3>
\beq
\label{ST}
S_m\^n(u;\hbar\>)\,=\,\sum_{k=0}^m\,\frac{(-1)^{m-k}}{(m-k)\<\>!}\;
T_k\^n(u;p,\hbar\>)\,\prod_{i=1}^{m-k}\,(p-m+i)\,.
\vv.3>
\eeq
\begin{rem}
Notice that \,$(1+x)^p=\sum_{m=0}^\infty\Tr_m\@p(A\^m)\,x^m$\>.
\end{rem}
\begin{rem}
Taking \,$p=m-1$ in either Proposition~\ref{TS} or formula~\Ref{ST} yields
\;$S_m\^n(u;\hbar)=T_m\^n(u;m-1;\hbar)$\>.
\end{rem}

\begin{lem}
\label{TNN}
Let \,$m\in\Z_{\ge\>n}$\,. Then
\vvn-.2>
\beq
\label{Tnmm}
T_m\^n(u;m)\,=\,\prod_{a=1}^n\,(u-z_a+\hbar\<\>)
\vv-.4>
\eeq
and \;$T_k\^n(u;m)=0$ \,for \,$k>m$.
\end{lem}
\begin{proof}
Since $m\ge n$, \>the homomorphism
\,$\varpi_n:\C[\Sg_n]\to\End\bigl((\C^m)^{\ox n}\bigr)$ is injective.
Hence, it suffices to compute the image of \,$T_k\^n(u;m)$ \,under
\,$\varpi_n$.
This can be done using Proposition~\ref{trace}. Since \,$\varpi_m(A\^k)=0$
\,\,for \,$k>m$\>, verifying that \;$T_k\^n(u;m)=0$ \,for \,$k>m$ \,is
straightforward. This equality also \,follows from Proposition~\ref{TS},
because \,$S_k\^n(u)=0$ \,for \,$k>n$.

\vsk.2>
To prove formula~\Ref{Tnmm}, we observe that \,$\varpi_m(A\^m)$ \,is
a rank-\<\>one projector and use Proposition~\ref{trace} to get
\vvn-.3>
\beq
\label{Trmm}
\Tr_m\@m\Bigl(\thi_m(A\^m)\>\tsum_{i=1}^m\si_{a,n+i}\;X\Bigr)\,=\,
\Tr_m\@m\Bigl(\thi_m(A\^m)\,X\Bigr)
\vv.2>
\eeq
for any \,$a=1\lc n$, \,and any \>$X\<\in\Sg_{n+m}$\>.
\,Together with formula~\Ref{Tnm}, this yields the claim.
\end{proof}

Taking \,$p=m$ \,in Proposition~\ref{TS}, we have
\beq
\label{Snn}
\sum_{k=0}^n\,S_k\^n(u;\hbar\>)\,=\,\prod_{a=1}^n\,(u-z_a+\hbar\<\>)\,.
\eeq

\vsk-.1>
Write
\vvn-.3>
\beq
\label{Snmi}
S_m\^n(u;\hbar\>)\,=\,\sum_{i=0}^{n-m}\>S_{\mi}\^n\,u^{n-m-i}\,.
\vvn-.3>
\eeq
Notice that by formula~\Ref{Sk},
\beq
\label{Sm0}
S_{\%m0}\^n\,=\,\hbar^m\>\Pho_{\%m0}\^n\,.
\eeq

\begin{prop}
\label{BS}
The subalgebra \,$\Bc^\ss_\nh\pz$ \,is generated by the elements
\,$S_{\mi}\^n$ \,for \;$m=1\lc n-1$, \,$i=1\lc n-m$.
\end{prop}
\begin{proof}
By formula~\Ref{ST}, the subalgebra \,$\Bc^\ss_\nh\pz$ \,contains
the elements \,$S_{\mi}\^n$ \,for all \,$m$ \,and \,$i$\>. Recall that
\;$S_0\^n(u;\hbar\>)=\prod_{a=1}^n(u-z_a)$. By Proposition~\ref{TS},
the subalgebra \,$\Bc^\ss_\nh\pz$ \,is generated by the elements
\vvn.1>
\,$S_{\mi}\^n$ \,for \;$m=1\lc n-1$\>, \,$i=0\lc n-m$\>,
\,and by formula~\Ref{Snn},
\;$S_{\%m0}\^n=\<\>c_m-\sum_{k=1}^{m-1}\,S_{k,\<\>m-k}\^n$\>, \>where
\vvn.1>
\be
\sum_{m=1}^{n-1}\<\>c_m\>u^{n-m}\<\>=\,
\prod_{a=1}^n\,(u-z_a\<+\hbar)\>-\prod_{a=1}^n\,(u-z_a)\,.
\vvn->
\ee
\end{proof}

\begin{cor}
The subalgebra \,$\Bc^\ss_\nh\pz$ \,contains the elements
\,$T_{\mi}\^n$\>, see~\Ref{Tmi}, \;for all \,$m$ and \,$i$\>.
\end{cor}
\begin{proof}
The claim follows from Proposition~\ref{TS}.
\end{proof}

\begin{lem}
\label{Bhsz}
We have \;$\Bc^\ss_{n,s\hbar}(sz_1\lc sz_n)=\Bc^\ss_\nh\pz$
\vvn.16>
\,for any \,$s\ne 0$, \;and\\
\;$\Bc^\ss_\nh(z_1+s\lc z_n+s)= \Bc^\ss_\nh\pz$ \,for any \,$s$.
\end{lem}
\begin{proof}
Formula~\Ref{Tnm} yields
\;$T_m\^n(su;p;s\hbar\<\>;sz_1\lc sz_n)=s^n\>T_m\^n(u;p;\hbar\<\>;\zn)$.
\vvn.16>
Hence, \;$T_{\mi}\^n(\<\>p;s\hbar\<\>;sz_1\lc sz_n)=
s^i\>T_{\mi}\^n(\<\>p;\hbar\<\>;\zn)$, \,which proves the first claim.
Similarly, the second claim follows from the equality
\;$T_m\^n(u+s;p;\hbar\<\>;z_1+s\lc z_n+s)=T_m\^n(u;p;\hbar\<\>;\zn)$.
\end{proof}

Let \,$\gm\<\in\Sg_n$ \,be given by \;$\gm(i)=i+1$ \,for \,$i=1\lc n-1$\>,
\,and \,$\gm(n)=1$\>.

\begin{prop}
\label{cycle}
We have \;$\gm\,T_m\^n(u;p;\hbar\<\>;z_2\lc z_n,z_1)\,\gm^{-1}=\>
T_m\^n(u;p;\hbar\<\>;\zn)$\>.
\end{prop}
\begin{proof}
By formula~\Ref{Tnm} and Lemma~\ref{TrXY}, we have
\vvn.4>
\begin{align*}
& T_m\^n(u;p;\hbar\<\>;\zn)\,={}
\\[4pt]
&\!{}=\,\Tr_m\@p\Bigl(\bigl(u-z_1+\hbar\tsum_{i=1}^m\si_{1,n+i}\bigr)\,
\thi_m(A\^m)\>\bigl(u-z_n+\hbar\tsum_{i=1}^m\si_{n,n+i}\bigr)\dots
\bigl(u-z_2+\hbar\tsum_{i=1}^m\si_{2,n+i}\bigr)\<\Bigr)\,={}
\\[3pt]
&\!{}=\,\Tr_m\@p\Bigl(\thi_m(A\^m)\>
\bigl(u-z_1+\hbar\tsum_{i=1}^m\si_{1,n+i}\bigr)\>
\bigl(u-z_n+\hbar\tsum_{i=1}^m\si_{n,n+i}\bigr)\dots
\bigl(u-z_2+\hbar\tsum_{i=1}^m\si_{2,n+i}\bigr)\<\Bigr)\,={}
\\[4pt]
&\!{}=\,\gm\,T_m\^n(u;p;\hbar\<\>;z_2\lc z_n,z_1)\,\gm^{-1}\,.
\end{align*}
In the second equality we use that \,$\tsum_{i=1}^m\si_{1,n+i}\;\thi_m(A\^m)
\>=\,\thi_m(A\^m)\,\tsum_{i=1}^m\si_{1,n+i}$\>.
\end{proof}

\begin{cor}
We have
\;$\Bc^\ss_\nh(z_2\lc z_n,z_1)\>=\>\gm^{-1}\,\Bc^\ss_\nh\pz\,\gm$\>.
\end{cor}

\begin{prop}
\label{Tzz}
For any \,$a=1\lc n$, \>we have
\vvn.3>
\begin{align*}
\bigl((z_a-z_{a+1})\>\si_{a,a+1}+\hbar\>)\>
T_m\^n(u;p\<\>;\hbar\<\>;\zn) &{}\,={}
\\[4pt]
{}=\,T_m\^n(u;p\<\>;\hbar\<\>;z_1\lc z_{a+1},z_a\lc z_n) &
\>\bigl((z_a-z_{a+1})\>\si_{a,a+1}+\hbar\>)\,.
\\[-8pt]
\end{align*}
\end{prop}

\begin{proof}
The claim follows from the identity
\begin{align*}
\bigl(( & z_a-z_{a+1})\>\si_{a,a+1}+\hbar\>)\>
\bigl(u-z_{a+1}+\hbar\tsum_{i=1}^m\si_{a+1,n+i}\bigr)\<\Bigr)\>
\bigl(u-z_a+\hbar\tsum_{i=1}^m\si_{a,n+i}\bigr)\<\Bigr)\,={}
\\
& \!\<{}=\,\bigl(u-z_a+\hbar\tsum_{i=1}^m\si_{a+1,n+i}\bigr)\<\Bigr)\>
\bigl(u-z_{a+1}+\hbar\tsum_{i=1}^m\si_{a,n+i}\bigr)\<\Bigr)\>
\bigl((z_a-z_{a+1})\>\si_{a,a+1}+\hbar\>)\,.
\\[-40pt]
\end{align*}
\vv1.1>
\end{proof}

\begin{cor}
\label{Bzzh}
For any \,$a=1\lc n$, \>we have
\vvn.3>
\begin{align*}
\bigl((z_a-z_{a+1})\>\si_{a,a+1}+\hbar\>)\>\Bc^\ss_\nh\pz &{}\,={}
\\[4pt]
{}=\,\Bc^\ss_\nh(z_1\lc z_{a+1},z_a\lc z_n) &
\>\bigl((z_a-z_{a+1})\>\si_{a,a+1}+\hbar\>)\,.
\\[-12pt]
\end{align*}
\end{cor}

\begin{cor}
\label{conjprop}
If \,$z_a\ne z_{a+1}\pm\hbar$ \,then the subalgebras
\;$\Bc^\ss_\nh\pz$ \,and \;$\Bc^\ss_\nh(z_1\lc z_{a+1},z_a\lc z_n)$
\,are conjugate in \,$\C[\Sg_n]$\>, cf.~\Ref{Bzz}.
\end{cor}
\begin{proof}
Since \,$\bigl((z_a-z_{a+1})\>\si_{a,a+1}+\hbar\>)\>
\bigl((z_a-z_{a+1})\>\si_{a,a+1}-\hbar\>)=(z_a-z_{a+1})^2-\hbar^2\ne 0$,
\vvn.16>
\;the claim follows from Corollary~\ref{Bzzh}.
\end{proof}

\begin{prop}
The subalgebra \,$\Bc^\ss_\nh\pz$ contains the center of
\,$\C[\Sg_n]$.
\end{prop}
\begin{proof}
The statement follows from formula~\Ref{Sm0} and Proposition~\ref{center}.
\end{proof}

\begin{lem}
\label{tlem}
We have
\vvn.3>
\begin{align*}
\thi_m(A\^m)\>\bigl(u-(m-1)\>\hbar+\hbar\>\si_{a,n+1}\bigr)\,\dots\,
\bigl(u+\hbar\>\si_{a,n+1}\bigr)& {}\,={}
\\[2pt]
\thi_m(A\^m)\>\bigl(u+\hbar\tsum_{i=1}^m\si_{a,n+i}\bigr) &
\,\prod_{i=1}^{m-1}\,(u-i\<\>\hbar)\,.
\end{align*}
\end{lem}
\begin{proof}
Both sides of the formula are polynomials in \>$u$ \>of degree \>$m$
with the matching coefficients for \>$u^m$ and \>$u^{m-1}$. In addition,
the product \,$\prod_{i=1}^{m-1}(u-i\<\>\hbar)$ \,divides the left-hand side
of the formula due to the identity
\vvn.3>
\be
(1-\si_{i,\<\>j})\>(u-\hbar+\hbar\>\si_{a,\<\>i})\>(u+\hbar\>\si_{a,j})\,=\,
(u-\hbar\>)\>(1-\si_{i,\<\>j})\>(u+\hbar\>\si_{a,\<\>i}+\hbar\>\si_{a,j})\,,
\vv.4>
\ee
which completes the proof.
\end{proof}

\begin{thm}
\label{Bhcomm}
The subalgebra \,$\Bc^\ss_\nh\pz$ is commutative.
\end{thm}
\begin{proof}
Let \,$P(u)=\prod_{a=1}^n\,(u-z_a)$\>.
\,Since \,$\thi_m(A\^m)$ \,commutes with
\,$\bigl(u+\hbar\tsum_{i=1}^m\si_{a,n+i}\bigr)$\>,
formula~\Ref{Tnm} can be written as
\begin{align*}
T_m\^n(u;p\<\>;\hbar\>)\>=\>P(u)\;\Tr_m\@p\biggl(\thi_m(A\^m)\>
\Bigl(1+\frac{\hbar\>\si_{n,n+1}}{u-z_n-(m-1)\>\hbar}\>\Bigr)\,\dots\,
\Bigl(1+\frac{\hbar\>\si_{n,n+m}}{u-z_n}\Bigr)&
\\[4pt]
{}\dots\;\Bigl(1+\frac{\hbar\>\si_{1,n+1}}{u-z_1-(m-1)\>\hbar}\>\Bigr)\,\dots\,
\Bigl(1+\frac{\hbar\>\si_{1,n+m}}{u-z_1}\Bigr)&{}\biggr)\,,\kern-.5em
\\[-14pt]
\end{align*}
see Lemma~\ref{tlem}. Hence, for \>$N\!\in\Z_{>0}$\>, \,the image of
\vvn.2>
\;$T_m\^n(\hbar\<\>u;N;\hbar\>)/P(\hbar\<\>u)$ \>under the homo\-morphism
\;$\varpi_n:\C[\Sg_n]\to\End\bigl((\C^N)^{\ox n}\bigr)$ \,coincides with
the action of the \,$m$-th transfer-mat\-rix \,$\Te_m(u)$, see~\Ref{Tem}, on
the \,$\Yn$-module \,$V^\yv(z_1/\hbar)\ox\dots\ox V^\yv(z_n/\hbar)$. Therefore,
\vvn.3>
\be
\varpi_n\bigl(\Bc^\ss_\nh\pz\bigr)\,=\,\Bc^\yy_\nnla(z_1/\hbar\lc z_n/\hbar)\,,
\vv.3>
\ee
cf.~Theorem~\ref{SW}. Here \,$\Bc^\yy_\nnla(z_1/\hbar\lc z_n/\hbar)$ is defined
in Section~\ref{A3}. Since the homomorphism \,$\varpi_n$ \,is injective for
\,$N\ge n$, the theorem follows from Theorem~\ref{BY}.
\end{proof}

\begin{prop}
\label{BBh}
Let \,$\zn$ be distinct. Then the subalgebra
\,$\Bc^\ss_\nh\pz$ \,tends to \,$\Bc^\ss_n\pz$
\,as \,$\hbar\to 0$.
\end{prop}
\begin{proof}
By formula~\Ref{Sk} and Proposition~\ref{BS}, the limit of \,$\Bc^\ss_\nh\pz$
contains \,$\Bc^\ss_n\pz$. Since \,$\Bc^\ss_n\pz$ is a maximal commutative
subalgebra of \,$\C[\Sg_n]$ \>for distinct \,$\zn$, and \,$\Bc^\ss_\nh\pz$
\,is commutative for any $\zn$, the limit of \,$\Bc^\ss_\nh\pz$ coincides with
\,$\Bc^\ss_n\pz$.
\end{proof}

\begin{prop}
\label{GZh}
The Bethe subalgebra \,$\Bc^\ss_\nh\pz$ \,tends to the Gelfand\<\>-Zetlin
subal\-gebra \;$\Gc_n\?$ \,as \>$\hbar/(z_1-z_2)\to 0$ \;and
\,$(z_{a-1}-z_a)/(z_a-z_{a+1})\to0$ \;for \,$a=2\lc n-1$\>.
\end{prop}
\begin{proof}
Without loss of generality we can assume that $z_1=0$, see Lemma~\ref{Bhsz},
so we have \,$h/\<z_2\to0$ and \,$z_a\<\>/\<z_{a+1}\to0$
\;for all \,$a=2\lc n-1$\>. In this limit the element
\,$(-1)^j\>S_{\ij}\^n\,z_{n-j+1}^{-1}\dots z_n^{-1}$ \>tends to
$\pi_{1\lc n-j}\^n\bigl(\Pho_{\%i-j0}\^{n-j\<\>}\bigr)$, see formulae~\Ref{Sk},
\vvn.1>
\Ref{Snmi}, \Ref{Phoi0}. Hence, the limit of \,$\Bc^\ss_\nh\pz$ contains
\vvn.08>
\,$\Gc_n$. Since \,$\Gc_n$ is a maximal commutative subalgebra of
\,$\C[\Sg_n]$, and \,$\Bc^\ss_\nh\pz$ \,is commutative for any $\zn$,
\vvn.1>
the limit of \,$\Bc^\ss_\nh\pz$ coincides with \,$\Gc_n$.
\end{proof}

\begin{rem}
Recall that \,$\Zc_m$ \,is the center of \,$\C[\Sg_m]$\>.
Similarly to~\Ref{Tnm}, for any \>$X\<\in\Zc_m$, \>set
\vvn.2>
\be
T_{m,X}\^n(u;p\<\>;\hbar\>)\,=\,\Tr_m\@p\Bigl(\thi_m(X)\>
\bigl(u-z_n+\hbar\tsum_{i=1}^m\si_{n,n+i}\bigr)\dots
\bigl(u-z_1+\hbar\tsum_{i=1}^m\si_{1,n+i}\bigr)\<\Bigr)\,,
\ee
$T_{m,X}\^n(u;p\<\>;\hbar\>)=
\sum_{i=0}^n\,T_{m,X,\>i}\^n(p\<\>;\hbar\>)\>u^{n-i}$.
\,Then \,$T_{m,X,\>i}\^n(p\<\>;\hbar\>)\in\Bc^\ss_\nh\pz$ \,for all
\,$X\<\in\Zc_m$ and \,$i=1\lc n$\>. \,However, the constructed linear map
\,$\Zc_m\to\Bc^\ss_\nh\ox\C[u]$\,, \ $X\mapsto T_{m,X}\^n(u;p\<\>;\hbar\>)$\>,
\,is not a homomorphism of algebras.
\end{rem}

Let \,$\rhe\<\in\Sg_n$ \,be given by \;$\rhe(i)=n-i+1$ \;for \,$i=1\lc n$\>.

\begin{lem}
\label{Tnsi}
Let \>$N\in\Z_{\ge n}$. \,Then
\vvn.2>
\be
\rhe\,T_m\^n(u;N;\hbar\<\>;z_n\lc z_1)\,\rhe\,=\,
(-1)^n\,T_{N\?-m}\^n(-u-\hbar\<\>;N;\hbar\<\>;-z_1\lc-z_n)
\vv.3>
\ee
for all \,$m=0\lc n$.
\end{lem}
\begin{proof}
By Lemmas~\ref{TrA} and~\ref{TrXY}, \,formula~\Ref{Tnm} implies that
\vvn.2>
\begin{align}
\label{Tnk}
& T_k\^n(u;p,\hbar\>)\,\prod_{i=1}^{m-k}\,\frac{p-m+i}{m+1-i}\,={}
\\[4pt]
& \!{}=\,\Tr_m\@p\Bigl(\thi_m(A\^m)\>
\bigl(u-z_n+\hbar\tsum_{i\in K}\si_{n,n+i}\bigr)\dots
\bigl(u-z_1+\hbar\tsum_{i\in K}\si_{1,n+i}\bigr)\<\Bigr)\,,
\notag
\\[-16pt]
\notag
\end{align}
where \,$K$ \>is any \>$k$-element subset of \,$\{1\lc m\}$\>. Then we have
\vvn.5>
\begin{alignat*}3
& \hphantom{\,{}-m)\<\>!}\llap{$\dsty\frac{m\<\>!\,(N-m)\<\>!}{N\<\>!}$}
\rlap{$\;\rhe\,T_m\^n(u;N;\hbar\<\>;z_n\lc z_1)\,\rhe\,={}$}
\\[4pt]
\hphantom{m\<\>!\,(N} & {}=\,\Tr_N\@N\Bigl(&{}\>\thi_N(A\^N)\>
& \bigl(u-z_1+\hbar\tsum_{i=1}^m\si_{1,n+i}\bigr)\,
\dots\,\bigl(u-z_n+\hbar\tsum_{i=1}^m\si_{n,n+i}\bigr)\<\Bigr)\,={} &&
\\[4pt]
& {}=\,\Tr_N\@N\Bigl(&{}\>\thi_N(A\^N)\>
& \bigl(u-z_1+\hbar\tsum_{i=1}^N\si_{1,n+i}-
\hbar\!\tsum_{i=m+1}^N\si_{1,n+i}\bigr)\times{}
\\[2pt]
&&{}\times{}& \bigl(u-z_2+\hbar\tsum_{i=1}^m\si_{2,n+i}\bigr)\,
\dots\,\bigl(u-z_n+\hbar\tsum_{i=1}^m\si_{n,n+i}\bigr)\<\Bigr)\,={}
\\[3pt]
\noalign{\penalty-500}
& {}=\,\Tr_N\@N\Bigl(&{}\>\thi_N(A\^N)\>
& \bigl(u-z_1+\hbar-\hbar\!\tsum_{i=m+1}^N\si_{1,n+i}\bigr)\times{}
\\[1pt]
&&{}\times{}& \bigl(u-z_2+\hbar\tsum_{i=1}^m\si_{2,n+i}\bigr)\,
\dots\,\bigl(u-z_n+\hbar\tsum_{i=1}^m\si_{n,n+i}\bigr)\<\Bigr)\,={}
\\[5pt]
& {}=\,\Tr_N\@N\Bigl(&{}\>\thi_N(A\^N)\>
& \bigl(u-z_2+\hbar\tsum_{i=1}^m\si_{2,n+i}\bigr)\,
\dots\,\bigl(u-z_n+\hbar\tsum_{i=1}^m\si_{n,n+i}\bigr)\times{}
\\[2pt]
&&&& \llap{${}\times\bigl(u-z_1+\hbar-
\hbar\!\tsum_{i=m+1}^N\si_{1,n+i}\bigr)\<\Bigr)\,={}$} &
\\[-14pt]
\end{alignat*}
where for the third equality, we use formula~\Ref{Trmm}.
Repeating similar transformations several times, we arrive at
\vvn.2>
\begin{align*}
&\!{}=\,\Tr_N\@N\Bigl(\thi_N(A\^N)\>
\,\bigl(u-z_n+\hbar-\hbar\!\tsum_{i=m+1}^N\si_{n,n+i}\bigr)\,
\dots\,\bigl(u-z_1+\hbar-\hbar\!\tsum_{i=m+1}^N\si_{1,n+i}\bigr)\<\Bigr)\,={}
\\[3pt]
&\!{}=\,(-1)^n\,\frac{m\<\>!\,(N-m)\<\>!}{N\<\>!}\;
T_{N-m}\^n(-u-\hbar\<\>;N;\hbar\<\>;-z_1\lc-z_n)\,,
\\[-12pt]
\end{align*}
the last equality following from~\Ref{Tnk}. The proposition is proved.
\end{proof}

\begin{cor}
\label{B-}
We have
\;$\Bc^\ss_\nh(-z_1\lc -z_n)\>=\>\rhe\,\Bc^\ss_\nh(z_n\lc z_1)\,\rhe$\>.
\end{cor}

Recall that \,$^\dag$ and \,$^*$ are the linear and semilinear antiinvolutions
on $\C[\Sg_n]$ such that \,$\si^\dag\<=\si^*\<=\si^{-1}$ for any $\si\in\Sg_n$.

\begin{prop}
\label{T+*}
Let \>$N\in\Z_{\ge n}$. \,Then
\vvn.2>
\begin{align}
\bigl(T_m\^n(u;N;\hbar\<\>;\zn)\bigr)^{\<\dag}\> &{}=\,
(-1)^n\,T_{N\?-m}\^n(-u-\hbar\<\>;N;\hbar\<\>;-z_1\lc-z_n)\,,
\notag
\\[4pt]
\bigl(T_m\^n(u;N;\hbar\<\>;\zn)\bigr)^{\<*}\<\> &{}=\,
(-1)^n\,T_{N-m}\^n(-\bar u-\bar\hbar\<\>;N;\bar\hbar\<\>;-\bar z_1\lc-\bar z_n)\,,
\label{T*}
\\[-12pt]
\notag
\end{align}
for all \,$m=0\lc n$\>. \,Here \,$\bar u\<\>,\bar\hbar\<\>,\bar z_1\lc\bar z_n$
\>are the complex conjugates of \,$u\<\>,\hbar\<\>,\zn$\>.
\end{prop}
\begin{proof}
It is easy to see that
\vvn-.3>
\beq
\label{Tr+}
\Tr_k\@p(X^\dag)\>=\,\bigl(\Tr_k\@p(X)\bigr)^{\<\dag}
\vv.2>
\eeq
for any \>$X\<\in\Sg_{n+k}$\>, where we use \;$^\dag$ \,for both
antiinvolutions of \,$\C[\Sg_{n+k}]$ \,and \,$\C[\Sg_n]$.
\vvn.3>
Since
\beq
\label{T+}
\bigl(T_m\^n(u;p,\hbar\<\>;\zn)\bigr)^\dag\>=\,
\rhe\,T_m\^n(u;p,\hbar\<\>;z_n\lc z_1)\,\rhe\,,
\vv.5>
\eeq
by~\Ref{Tnm} and~\Ref{Tr+}, the first equality follows from Lemma~\ref{Tnsi}.
The second equality is then straightforward.
\end{proof}

\begin{cor}
\label{Bh+*}
We have
\vvn.7>\\
\leftline{\hfil$
\bigl(\Bc^\ss_\nh\pz\bigr)^{\<\dag}=\>\Bc^\ss_\nh(-z_1\lc-z_n)$\quad
and\quad
$\bigl(\Bc^\ss_\nh\pz)^{\<*}=\>
\Bc^\ss_\nhb(-\bar z_1\lc-\bar z_n)$\>.}
\end{cor}

\section{Further properties of the Bethe subalgebras
\,$\Bc^\ss_\nh\pz$}
\label{morexxx}

For a partition \,$\bs\la$ \>of \>$n$, recall that $M_{\bs\la}$ denotes
the irreducible $\Sg_n$-module corresponding to $\bs\la$, and
$\chip_{\bs\la}\<\in\C[\Sg_n]$ \,--- \>the respective central idempotent.
Set \,$\Bc^\ss_\nhla\pz=\chip_{\bs\la}\Bc^\ss_\nh\pz$.
The algebra \>$\Bc^\ss_\nhla\pz$ is isomorphic to the image of
\>$\Bc^\ss_\nh\pz$ in \,$\End(M_{\bs\la})$ by the canonical projection.
\,Clearly, \,$\Bc^\ss_\nh\pz=\tbigoplus_\lapns\Bc^\ss_\nhla\pz$.

The next lemma is straightforward.

\begin{lem}
The assignment
\vvn-.1>
\beq
\label{sidot}
\si_{\ii+1}:q(y_1\lc y_n)\,\mapsto\,
q(y_1\lc y_{i+1},y_i\lc y_n)\,\frac{y_i-y_{i+1}+\hbar}{y_i-y_{i+1}}\>-\>
\frac{\hbar\,q(y_1\lc y_n)}{y_i-y_{i+1}}
\vv.3>
\eeq
for \,$i=1\lc n-1$\>, \,defines an action of \;$\Sg_n\<$ on
\,$\C[y_1\lc y_n]$\>.
\qed
\end{lem}

Denote by \,$\C[y_1\lc y_n]_\hbar$ \>the obtained \>$\Sg_n$-module. Recall
that \,$\deg\>q$ \,denotes the homogeneous degree of \,$q(y_1\lc y_n)$\>.
We extend
\vvn.16>
the degree \,to \,$M_{\bs\la}\ox\alb\C[y_1\lc y_n]_\hbar$ \>trivially on
the first factor. Notice that action~\Ref{sidot} does not increase the degree,
\vvn.1>
and the leading part of~\Ref{sidot} acts by permuting the variables.
\vvn.16>
Therefore, for every \,$w\in(M_{\bs\la}\ox\alb\C[y_1\lc y_n]_\hbar)^{\Sg_n}$
we have \,$\deg\>w\ge\sum_{i=1}^n(i-1)\>\la_i$\>, \,and the component of
\vvn.16>
\,$(M_{\bs\la}\ox\C[y_1\lc y_n]_\hbar)^{\Sg_n}$ of degree
$\sum_{i=1}^n(i-1)\>\la_i$ \>is one-dimensional, see Section~\ref{more}.
Let \;$w_\hla$ \>be
a nonzero element of \,$(M_{\bs\la}\ox\alb\C[y_1\lc y_n]_\hbar)^{\Sg_n}$
of degree $\sum_{i=1}^n(i-1)\>\la_i$.

\vsk.3>
For a positive integer \,$m$ \,and a partition \,$\bs\la$ \,with at most \>$m$
\>parts, consider indeterminates \>$f_{\ij}$ with \,$i=1\lc m$ \,and
\,$j=1\lc\la_i+m-i$, \;$j\ne\la_i-\la_s-i+s$ \,for \,$s=i+1\lc m$.
Given in addition a collection of complex numbers \,$\bs a=(a_1\lc a_n)$,
define the algebra $\Oc_\mlba$ as the quotient of
\,$\C[f_{1,1}\lc f_{m,\la_m}]$ \,by relations~\Ref{WOlah} described below.

\vsk.2>
Consider $f_{\ij}$ as coefficients of polynomials in one variable,
\be
f_i(u)\,=\,u^{\la_i+m-i}+\sum_{j=1}^{\la_i+m-i} f_{\ij}\,u^{\la_i+m-i-j}\,,
\vv-.2>
\ee
assuming that \,$f_{i,\la_i-\la_s-i+s}=0$. The last condition
\vvn.1>
means that the coefficient of \,$u^{\la_s+m-s}$ in \>$f_i(u)$ equals zero.
The defining relations for $\Oc_\mhlba$ can be written as an equality of
two polynomials in \,$u$:
\vvn-.3>
\beq
\label{WOlah}
\Wrh[f_1(u)\lc f_m(u)]\,=\prod_{1\le i<j\le m}\!(\la_i-\la_j-i+j)\;
\Bigl(u^n+\sum_{s=1}^n\,(-1)^s\,a_su^{n-s}\Bigr)\,,
\vv.2>
\eeq
where \;$\Wrh[f_1(u)\lc f_m(u)]=
\det\<\>\bigl(f_j(u-\hbar\>(i-1))\bigr)_{i,j=1\lc m}$
\;is the Casorati determinant (aka~the discrete Wronskian).
Notice that
\vvn.2>
\beq
\label{WOlahd}
\Wrh[f_1(u)\lc f_m(u)]\,=\,(-1)^{m(m-1)/2}\,
\det\<\>\bigl((1-e^{-\hbar\>\der_u})^{i-1}f_j(u)\bigr)_{i,j=1\lc m}\,.
\kern-2em
\vv.2>
\eeq

\begin{lem}
\label{Okmh}
Let \,$\bs\la$ \,be a partition with at most \>$m$ \>parts, and \,$k\ge m$.
Then the algebras \,$\Oc_\khlba$ \>and \,$\Oc_\mhlba$ \>are isomorphic.
\end{lem}
\begin{proof}
The proof is similar to the proof of Lemma~\ref{Okm}.
The counterpart of formula~\Ref{fi} is
\vvn.2>
\be
f_i^{\{m\}}(u)\,\mapsto\,
(1-e^{-\hbar\>\der_u})^{k-m}f_i^{\{k\}}(u)
\,\prod_{s=m+1}^k\,\frac1{\hbar\>(\la_i+s-i)}\,,
\vv-.3>
\ee
so that by \Ref{WOlahd},
\vvn-.3>
\be
\Wrh[f_1^{\{m\}}(u)\lc f_m^{\{m\}}(u)]\,\mapsto\,
\Wrh[f_1^{\{k\}}(u)\lc f_k^{\{k\}}(u)]\,
\prod_{s=0}^{k-m-1}\frac1{s\<\>!}
\,\;\prod_{i=1}^m\,\frac{(\la_i+m-i)\<\>!}{(\la_i+k-i)\<\>!}\;.
\vvn->
\ee
\end{proof}

Let \,$\Oc_\hlba=\Oc_\nhlba$. Set
\vvn.3>
\beq
\label{Fuv}
F(u,v)\,=\,e^{(n\hbar\>-u)x}\Wrh[f_1(u)\lc f_n(u),e^{ux}\>]\,,\qquad
v=e^{\<\>\hbar\>x}\,.
\vv.3>
\eeq
It is a polynomial in \,$u,v$ \,with coefficients in \,$\Oc_\hlba$\>:
\vvn.1>
\be
F_\hlba(u,v)\,=\,\sum_{i=0}^n\,
\sum_{m=0}^n\,(-1)^i\,F_\hlbami\,u^{\<\>n-i}\>v^{\<\>n-m}\,.
\vv.3>
\ee

Recall that we identify elements of \,${M_{\bs\la}\ox\alb\C[y_1\lc y_n]}$
\,with \,$M_{\bs\la}$-valued polynomials in \>$y_1\lc y_n$\>.

\begin{thm}
\label{manyproph}
Let \,$z_i-z_j\ne\hbar$ \>for all \;$1\le j<i\le n$. Then
\begin{enumerate}
\item[i)]
The algebra \,$\Bc^\ss_\nh\pz$ is a maximal commutative subalgebra of
\,$\C[\Sg_n]$.
\vsk.2>
\item[ii)]
The map \,$\Bc^\ss_\nh\pz\to\tbigoplus_\lapns M_{\bs\la}$\>,
\vvn-.1>
\;$X\mapsto\tbigoplus_\lapns Xw_\hla\pz$\>, is an isomorphism
of the regular representation of \,$\Bc^\ss_\nh\pz$ on itself
\vvn.1>
and the \,$\Bc^\ss_\nh$-module $\tbigoplus_\lapns M_{\bs\la}$\>.
In particular, \,$\dim\Bc^\ss_\nh\pz=\sum_\lapns\dim M_{\bs\la}$\;.
\vsk.4>
\item[iii)]
The action of \,$\Bc^\ss_\nh\pz$ on \,$\tbigoplus_\lapns M_{\bs\la}$
is diagonalizable and has simple spectrum if one of the following assumptions
holds:
\vsk.1>
\begin{enumerate}
\item[a)]
$(z_i-z_j)/\hbar$ \>are real for all \;$i,j=1\lc n$, and
there exists an integer \,$m$ \>such that \,$(z_i-z_{i+1})/\hbar>1$
\vvn.1>
\,for all \;$i=1\lc n-1$, \;$i\ne m$\>;
\vsk.1>
\item[b)]
There exists \,$y\<\in\C$ such that \,$(z_1-y)/\hbar\>\lc(z_n-y)/\hbar$
\>are either real or form pairs of complex conjugated numbers, and
\,$\bigl|\>\Im\>((z_i-y)/\hbar\>)\bigr|<1/2$
\,for all \,$i=1\lc n$\>.
\vsk.1>
\item[c)]
$\,\zn$ are generic.
\end{enumerate}
\vsk.2>
\item[iv)]
The assignment \,$T_{\mi}\^n(n;\hbar)\mapsto F_\hlbami$, \,for \,$i,m=1\lc n$,
extends to an isomorphism of algebras \,$\Bc^\ss_\nhla\pz\to\Oc_\hlba$.
\vvn.16>
Here \,$T_{\mi}\^n(n;\hbar)$ \,are given by~\Ref{Tmi}, \,$\bs a=(a_1\lc a_n)$
\>and
\;$u^n+\sum_{s=1}^n\,(-1)^s\,a_su^{n-s}=\prod_{i=1}^n(u-z_i+\hbar\<\>)$\,.
\vsk.4>
\item[v)]
The algebra \,$\Bc^\ss_\nhla\pz$ \,is a Frobenius algebra.
\end{enumerate}
\end{thm}
\begin{proof}
The proof of Theorem~\ref{manyproph} is similar to that
of Theorem~\ref{manyprop}. Due to the decomposition
$\Bc^\ss_\nh\pz=\tbigoplus_\lapns\Bc^\ss_\nhla\pz$, \,it suffices to verify
the counterparts of the claims for the algebras $\Bc^\ss_\nhla\pz$ and
the $\Bc^\ss_\nhla\pz$-modules $M_{\bs\la}$. The required statements follow
from the properties of the algebras \,$\Bc^\yy_\nnla(\xn)$ \,studied in
\cite{MTV10}. Recall that the algebra \,$\Bc^\yy_\nnla(\xn)$ coincides with the
\vvn.1>
algebra \,$\Bc^{\<\>\bs q=\bs1}\bigl(V(\bs x)^{\textit{sing}}_{\bs\la}\bigr)$,
where \,$\bs x=(\xn)$, in \cite{MTV10}.
\end{proof}

\begin{conj}
The assertions \;{\rm i)}\>, \,{\rm iv)}\>, \,{\rm v)}
\>of Theorem~\ref{manyproph} \,hold for all \,$\zn$.
\end{conj}

The next statement, which is similar to Theorem~2.1 in~\cite{MTV5}, is a rather
unexpected byproduct of Theorem~\ref{manyproph}. It is discussed in more detail
in \cite{MTV9}.

\begin{thm}
\label{Shapiro}
Let polynomials \,$p_1(t)\lc p_N(t)$ \,be such that the polynomial
\vvn.3>
\beq
\label{Wt}
W(t)\,=\,\bigl(\sqrt{-1}\,\<\>\bigr)^{N(N-1)/2}
\det\<\>\bigl(\<\>p_i\bigl(t+(N+1-2\<\>j\<\>)\>\sqrt{-1}\,\<\>\bigr)
\bigr)_{i,j=1\lc N}
\vv.2>
\eeq
has real coefficients, and all roots of \,$W(t)$ \>lie in the strip
\,$|\Im t\,|\le 1$. Then the vector space \>$X\?\subset\C[\<\>t\<\>]$\>,
\,spanned by the polynomials \,$p_1(t)\lc p_N(t)$, \,has a basis consisting
of polynomials with real coefficients.
\end{thm}
\begin{proof}
The statement follows from assertions~(ii), (iii\<\>,\,b), and~(iv) of
Theorem~\ref{manyproph}, formula~\Ref{T*}, and Proposition~\ref{Tzz},
\,provided \,$\hbar=2\>\sqrt{-1}$\;.
\end{proof}

\begin{rem}
If \,$p_1(t)\lc p_N(t)$ \,have real coefficients, then \,$W(t)$ \,has real
coefficients too.
\end{rem}

\begin{example}
Let \,$N=2$, \;$p_1(t)=t+a$\>, \;$p_2(t)=t^3+b\<\>t^2+c$\>, and \;$W(t)$ is
given by~\Ref{Wt}. Suppose \;$W(t)=4\<\>t^3+A\<\>t^2+B\<\>t$, \,with real
\,$A$ \>and \>$B$\>. Then the numbers \;$a,b,c$ \,are real if and only if
\,$B\le 4+A^2\?/\<\>12$\>, whereas Theorem~\ref{Shapiro} asserts that the
numbers \;$a,b,c$ \,are real provided \,$B\le 4+A^2\?/\<\>16$\>.
\end{example}

Consider the matrix \,$Z$ \,given by~\Ref{Z} and the matrix
\;$Q_\hbar$ \,with the entries
\vvn.2>
\beq
\label{Qh}
(Q_\hbar)_{ab}\,=\,\frac{\co_a}{z_a-z_b+\hbar}\;,
\vv.2>
\eeq
where \,$\co_1\lc\co_n$ \,are new variables.
Given a polynomial \,$q\in\C[u]$\>, \,set
\vvn.2>
\beq
\label{caz}
\co_a\,=\,q(z_a)\,\prod_{\fratop{b=1}{b\ne a}}^n\,\frac1{z_a-z_b}\;,
\qquad a=1\lc n\,,
\vvn-.8>
\eeq
and define
\vvn.1>
\beq
\label{Peh}
\Pe_\hbar(u,v;\zn;q)\,=\,
\det\<\>\bigl((u-Z)\>(v-Q_\hbar)-\hbar\>Q_\hbar\bigr)\,,
\vv.4>
\eeq
where \,$Q_\hbar$ \,is given by~\Ref{Qh}, \Ref{caz}.
Then \,$\Pe_\hbar(u,v;\zn;q)$ \,depends polynomially on the coefficients of
\,$q(u)$ \,and rationally on \,$\zn$ \,with possible poles at the hyperplanes
\,$z_a\<=z_b$ \,and \,$z_a\<+\hbar=z_b$ \,for \,$a\ne b$.

\begin{lem}
\label{Pehzz}
The expression \,$\Pe_\hbar(u,v;\zn;q)$ \,is regular at the hyperplanes
\,$z_a\<=z_b$\>.
\end{lem}
\begin{proof}
It follows from~\Ref{Peh} and~\Ref{caz} that \,$\Pe_\hbar(u,v;\zn;q)$
\,is invariant under permutations of \,$\zn$. So it suffices to show
that \,$\Pe_\hbar(u,v;\zn;q)$ \,is regular at the hyperplane
\,$z_1=z_2$. The order of the pole in question is at most two, and it is
straightforward to see that the coefficient of \,$(z_1-z_2)^{-2}$ actually
vanishes. Since \,$\Pe_\hbar(u,v;\zn;q)$ \,is invariant under the
exchange of \>$z_1$ and \>$z_2$, the coefficient of \,$(z_1-z_2)^{-1}$
vanishes too, which proves the lemma.
\end{proof}

Write
\vvn-.6>
\beq
\label{S1i}
S_1\^n(u;\hbar\>)\,=\,\hbar\>n\<\>u^{n-1}+\>
\sum_{i=1}^{n-1} \hbar^{\>i+1}\>S_{1,\<\>i}\^n\;u^{n-i-1}\,,
\vv.1>
\eeq
see~\Ref{S1}. \,Recall that
\;$T_1\^n(u;p\<\>;\hbar\>)=S_1\^n(u;\hbar\>)+p\>\prod_{a=1}^n(u-z_a)$\>,
\;see Proposition~\Ref{TS}.

\begin{thm}
\label{generateSh}
Let \,$z_a-z_b\ne\hbar$ \,for all \,$a,b=1\lc n$. Then the subalgebra
\,$\Bc^\ss_\nh\pz$ is generated by the elements
\,$S_{1,1}\^n\<\lc S_{1,n-1}\^n$. More precisely\>,
\vvn.1>
\beq
\label{TPeh}
T\^n(u,v;\hbar\>)\,=\,\Pe_\hbar(u,v;\zn;S_1\^n\>)\,,
\vv.2>
\eeq
where
\vvn-.3>
\beq
\label{Tuv}
T\^n(u,v;\hbar\>)\,=\,\sum_{m=0}^n\,(-1)^m\,T_m\^n(u;n;\hbar\>)\>v^{n-m}\,=\,
\sum_{k=0}^n\,(-1)^k\,S_k\^n(u;\hbar\>)\>(v-1)^{n-k}\,,
\vv-.3>
\eeq
\,see~\Ref{Tmi}, \Ref{Sk}.
\end{thm}
\begin{proof}
For distinct \,$\zn$\>, the proof is similar to the proof of Theorem~3.2
in~\cite{MTV7}, using the results from~\cite{MTV3}.
In this case, \,$\Bc^\ss_\nh\pz$ is generated by the elements
\,$K_a\^n\?=T_1\^n(z_a;\hbar\>)=S_1\^n(z_a;\hbar\>)$ 
\,for $a=1\lc n$, and hence, by the elements \,$S_{1,1}\^n\<\lc S_{1,n-1}\^n$,
because \,$S_1\^n(u;\hbar\>)$ \,is a polynomial in \>$u$ \>of degree \>$n-1$\>.
Since both sides of equality~\Ref{TPeh} are regular at the hyperplanes
\>$z_a\<=z_b$\>, \,see Lemma~\ref{Pehzz}, the statement follows.
\end{proof}

\begin{rem}
We have \;$T\^n(u,v;\hbar\>)=v^n\>\Th(u,\<-\<\>v^{-1};n;\hbar\>)=
(v-1)^n\>\Sh\bigl(u,(1-v)^{-1};\hbar\>\bigr)$\>,
\,see~\Ref{ThSh}.
\end{rem}

\begin{rem}
Notice that
\vvn.3>
\begin{align}
\label{Ka}
\Hh_a\^n\<\>={}\,
(z_a-z_{a-1}+\hbar\>\si_{a-1,\<\>a})\,\dots{} & (z_a-z_1+\hbar\>\si_{1,\<\>a})
\\[3pt]
{}\times\,(z_a-z_n+\hbar\>\si_{a,n})\dots{} &
(z_a-z_{a+1}+\hbar\>\si_{a,\<\>a+1})\,,
\notag
\\[-14pt]
\notag
\end{align}
cf.~\Ref{qKZ}. We call \,$K_1\^n\<\lc K_n\^n$ \>the \>\qKZ/ {\it elements\/}.
They commute with each other. Since
\vvn.1>
\be
\Hh_1\^n\<\dots\,\Hh_n\^n=\>\prod_{\satop{a,b=1}{a\ne b}}^n\,(z_a-z_b+\hbar)\,,
\ee
the \>\qKZ/ elements are invertible if \,$z_a-z_b\ne\hbar$
\,for all \,$a,b=1\lc n$.
\end{rem}

For $\zn$ such that $z_a-z_b\ne\hbar$ for all $a,b$, and a partition
$\bs\la=(\la_1\lc\la_n)$ of $n$, define the algebra \,$\He_\hla\pz$
as the quotient of \,$\C[\>q_1\lc q_n]$ \,by relations~\Ref{relHh},
\Ref{relHlh} described below. Consider a polynomial
\,$q(u)=\sum_{i=1}^nq_i\>u^{n-i}$, \,and write
\be
\Pe_\hbar(u,v;\zn;q)\,=\,
\sum_{i,j=0}^n\,\Pe_\hij(\zn;q)\,u^{n-j}\>(v-1)^{n-i}\,,
\ee
see~\Ref{Peh}. The defining relations for \,$\He_\hla\pz$ are
\vvn.3>
\beq
\label{relHh}
\Pe_\hij(\zn;q)\,=\,0\,,\qquad 0\le j<i\le n\,,
\eeq
and
\vvn-.5>
\beq
\label{relHlh}
\sum_{i=0}^n\,\Pe_\hii(\zn;q)\prod_{j=i+1}^n(t+j)\,=\,
\prod_{j=1}^n\,(t-\la_j+j)\,,
\vv.2>
\eeq
where \,$t$ is a formal variable, cf.~\Ref{relH}, \Ref{relHl}.

\begin{conj}
\label{BSnHh}
Let $z_a-z_b\ne\hbar$ for all $a,b=1\lc n$. Then the assignment
\,$q_1\mapsto\hbar\<\>\chip_{\bs\la}$ \,and
\,$q_i\mapsto \chip_{\bs\la}S_{1,i-1}\^n$ \,for \,$i=2\lc n$, \,defines
an isomorphism of algebras \;$\He_\hla\pz$ \,and \,$\Bc^\ss_\nhla\pz$\>.
\end{conj}

\section{Homogeneous Bethe subalgebra \,$\Aes_n$ of \,$\C[\Sg_n]$}
\label{homo}

Consider the subalgebra \,$\Aes_n=\Bc^\ss_\nh(z_1\lc z_1)\subset\C[\Sg_n]$.
By Lemma~\ref{Bhsz}, it does not depend on \,$\hbar$ \,and \,$z_1$.
We call \,$\Aes_n$ the {\it homogeneous Bethe subalgebra\/} of \,$\C[\Sg_n]$.
By Theorem~\ref{manyproph}, the subalgebra \,$\Aes_n$ is a maximal commutative
subalgebra of \,$\C[\Sg_n]$ \,of dimension \;$\sum_\lapns\dim M_{\bs\la}$\>.
We also have
\,$\bigl(\Aes_n\bigr)^{\<\dag}=\>\bigl(\Aes_n\bigr)^{\<*}=\>\Aes_n$\>,
\;see Corollary~\ref{Bh+*}.

\vsk.2>
Further on throughout this section we assume that \,$\hbar=1$ \,and
\,$z_1=\dots=z_n=0$.

\vsk.2>
Denote by \,$\Tht$ \,the set of \,$n\<\>$-dimensional subspaces
\,$U\subset\C[u]$ with the property: $U$ has a basis \,$p_1\lc p_n$ \,such that
\be
\Wrh[\>p_1(u)\lc p_n(u)]\,=\,(u+1)^n\>.
\vv.3>
\ee
For \,$U\<\in\Tht$ with a basis \,$p_1\lc p_n$ as above, set
\vvn.4>
\beq
\label{Fuvp}
F_U(u,v)\,=\,e^{(n-u)x}\Wrh[\>p_1(u)\lc p_n(u),e^{ux}\>]\,,\qquad
v=e^{\<\>x}\,.
\vv.3>
\eeq
Clearly, \,$F_U(u,v)$ does not depend on a choice of the basis.

\vsk.2>
For a partition \,$\bs\la$ \,of \>$n$, \,define the subset
\,$\Tht_{\bs\la}\subset\Tht$ \,as follows: \,$U\<\in\Tht_{\bs\la}$
\,if \,$U$ has a basis \,$p_1\lc p_n$ \,such that \;$\deg\>p_i=\la_i+n-i$\>,
\;$i=1\lc n$.

\begin{thm}
\label{lessprop}
The action of \;$\Aes_n$ on \,$\tbigoplus_\lapns M_{\bs\la}$ \,is
diagonalizable and has simple spectrum. The eigenvectors of $\Aes_n$ are in
a bijection with the elements of \,$\Tht$. If \,$\w_U$ is an eigenvector of
\;$\Aes_n$, corresponding to $U\<\in\Tht$, and \,$T\^n(u,v)$ \,is given
by~\,\Ref{Tuv} with \,$\hbar=1$, then
\vvn.2>
\be
T\^n(u,v)\,\w_U\,=\,F_U(u,v)\,\w_U\,,
\vv.2>
\ee
The vector \,$\w_U$ lies in the direct summand
\,$M_{\bs\la}\subset\tbigoplus_\lapns M_{\bs\la}$ \,if and only if
\,$U\<\in\Tht_{\bs\la}$\>.
\end{thm}
\begin{proof}
The theorem follows from items \;iii) \>and \;iv) \>of
Theorem~\ref{manyproph}\>.
\end{proof}

\begin{rem}
For \>$k\le n$, denote by \,$\Tht_k$ \,the set of \,$k\<\>$-dimensional
subspaces \,$U\subset\C[u]$ with the property: $U$ has a basis \,$p_1\lc p_k$
\,such that
\be
\Wrh[\>p_1(u)\lc p_k(u)]\,=\,(u+1)^n\>.
\vv.3>
\ee
For \,$U\<\in\Tht_k$ with a basis \,$p_1\lc p_k$ as above, set
\vvn.3>
\be
F_U(u,v)\,=\,e^{(k-u)x}\Wrh[\>p_1(u)\lc p_k(u),e^{ux}\>]\,,\qquad
v=e^{\<\>x}\,.
\vv.2>
\ee
Set \,$\Tht_{(k)}=\>\bigcup_{\lapns,\<\>\la_{k+1}=0}\Tht_{\bs\la}$\>.
For \,$U\<\in\Tht_{(k)}$ \,with a basis \,$p_1\lc p_n$, \,let
\vvn.16>
\,$\dl_k(U)\subset\C[u]$ \,be the subspace, spanned by the polynomials
\,$(1-e^{-\der_u})^{n-k}\<\>p_i(u)$, \;$i=1\lc n$\>. It is easy to see that
$\dl_k$ is a bijection from \,$\Tht_{(k)}$ \,to \,$\Tht_k$ \,and
\;$F_U(u,v)\>=\>(v-1)^{n-k}\>F_{\dl_k(U)}(u,v)$\,.
\end{rem}
Set
\vvn-.2>
\beq
\label{Gk}
G_k\^n\,=\sum_{1\le\>i_1<\<\cdots<\>i_k\le n\!\!}
\si_{i_1,\<\>i_2}\>\si_{i_2,\<\>i_3}\dots\<\>\si_{i_{k-1},\<\>i_k}\,,
\qquad k=2\lc n\,.
\eeq
By formula~\Ref{S1} with \,$\hbar=1$, we have
\vvn-.2>
\beq
\label{S1G}
S_1\^n(u)\,=\,n\<\>u^{n-1}+\sum_{k=2}^n G_k\^n\>u^{n-k}\,.
\vv-.4>
\eeq
Recall that \;$T_1\^n(u;p\<\>;\hbar\>)=p\>u^n+S_1\^n(u;\hbar\>)$.

\vsk.3>
Let \,$\Zh$ be the matrix with entries \,$\Zh_{ab}=\dl_{a,b-1}$\>,
\;$a,b=1\lc n$\>. Given a polynomial \,$q\in\C[u]$\>, consider the matrix
\;$\Qh$ \,with entries
\vvn.3>
\beq
\label{Qht}
\Qh_{ab}\,=\,\frac1{(n-a)\>!}\,
\frac{d^{\>n-\<\>a}}{du^{\<\>n-\<\>a}}\>
\Bigl(\>\frac{q(u)}{(u+1)^b}\,\Bigr)\Big|_{u=0}
\vv-.3>
\eeq
and define
\vvn.2>
\beq
\label{Peht}
\Peh(u,v;q)\,=\,
\det\<\>\bigl((u-\Zh\>)\>(v-\Qh\>)-\Qh\>\bigr)\,,
\vvn.2>
\eeq
cf.~\Ref{Peh}.
\,$\Peh(u,v;q)$ \,depends polynomially on the coefficients of \,$q(u)$\>.

\begin{lem}
\label{PePe}
Let \,$\hbar=1$ \,and \;$z_1=\dots=z_n=0$\>.
Then \;$\Pe_\hbar(u,v;\zn;q)\>=\>\Peh(u,v;q)$\>.
\end{lem}
\noindent
The proof is given in Section~\ref{A4}.

\begin{thm}
\label{generateG}
The subalgebra \,$\Aes_n$ is generated by the elements
\,$G_2\^n\<\lc G_n\^n$. More precisely\>,
\;$T\^n(u,v)\>=\>\Peh(u,v;S_1\^n\>)$\>,
\,where \;$T\^n(u,v)$ \,is given by~\,\Ref{Tuv} with \,$\hbar=1$\>.
\end{thm}
\begin{proof}
Since \,$G_k\^n=\>S_{1,k-1}\^n$\>, \,see~\Ref{S1i}, \Ref{S1G},
the statement follows from Theorem~\ref{generateSh} and Lemma~\ref{PePe}.
\end{proof}

Let \;$\gm_n=\si_{1,2}\,\si_{2,3}\>\dots\<\>\si_{n-1,n}=\>G_n\^n$. \,Then
\vvn-.1>
\beq
\label{gmn}
\gm_n^{-1}\>S_1\^n(u)\,=\,1+\sum_{k=1}^{n-2}\,\gm_n^{-1}\>G_{n-k}\^n\>u^k
+\gm_n^{-1}\>n\<\>u^{n-1}
\vvn-.2>
\eeq
\,and there exists the power series
\vvn-.3>
\beq
\label{Ik}
\log\>\bigl(\gm_n^{-1}\>S_1\^n(u)\bigr)\,=\,
\sum_{k=1}^\infty\,I_k\^n\>u^k
\vv-.6>
\eeq
with coefficients in \,$\C[\Sg_n]$\>.

\vsk.2>
Recall that \,$\pi_{1\lc k}\^n:\C[\Sg_k]\to\C[\Sg_n]$ \,is the embedding
\vvn-.1>
induced by the correspondence \,$i\>\mapsto i$.
For brevity, set \,$\pi_k=\pi_{1\lc k}\^n$\,.

\begin{thm}
\label{Inkt}
For every \;${k\<\in\Z_{>0}}$ \,there exists an element
\;${\tht_k\?\in\Sg_{k+1}}$ independent of \;$n$ \<\>such that
\vvn-.2>
\beq
\label{Ink}
I_k\^n\,=\,\sum_{m=0}^{n-1}\,\gm_n^{\<\>m}\,\pi_{k+1}(\tht_k)\,\gm_n^{-m}\>,
\qquad k=1\lc n-2\,.
\vv.2>
\eeq
For example, one can take \;$\tht_1=\si_{1,2}$\>,
\;$\tht_2=(\si_{2,3}\,\si_{1,2}-\si_{1,2}\,\si_{2,3}-1)/2$\>.
\end{thm}
\begin{proof}
\,The proof is by the same arguments as in the Appendix of~\cite{L}.

\vsk.2>
Given integers \,$0=i_0\<<i_1<\<\dots<\>i_k\<<i_{k+1}\<=n+1$ \,with \,$k<n$
\,define \,$\vec\si_{i_1\lc i_k}\!\in\Sg_n$ \,as follows.
Choose \,$m\in\{1\lc k\}$ \,such that \,$i_{m+1}\<-i_m>1$ \,and set
\vvn.4>
\be
\vec\si_{i_1\lc i_k}=\,\si_{i_m,i_m+1}\dots\<\>\si_{i_1,i_1+1}\;
\si_{i_k,i_k+1}\dots\<\>\si_{i_{m+1},i_{m+1}+1}\,,
\vv.4>
\ee
where \,$\si_{n,n+1}$ is understood as \,$\si_{1,n}$. It is easy to see that
\vvn.1>
\,$\vec\si_{i_1\lc i_k}$ does not depend on the choice of \>$m$.
By~\Ref{Gk}, we have \;$\gm_n^{-1}\>G_{n-k}\^n=
\sum_{1\le\>i_1<\<\dots<\>i_k\le n}\vec\si_{i_1\lc i_k}$\>.
Consider the series
\beq
\label{logn}
\log\>\biggl(1+\sum_{k=1}^{n-2}\,\sum_{1\le\>i_1<\<\dots<\>i_k\le n\!\!}
\vec\si_{i_1\lc i_k}\,u_{i_1}\?\dots u_{i_k}\biggr)\,=
\sum_{s_1\lc s_n=0\!}^\infty\phi_{s_1\lc\>s_n}\,u_1^{s_1}\?\dots\>u_n^{s_n}
\vv.3>
\eeq
in the variables \,$u_1\lc u_n$ with coefficients in \,$\C[\Sg_n]$.
\vvn.16>
It is easy to see that
\,$\phi_{s_2\lc\>s_n,\>s_1}=\gm_n\>\phi_{s_1\lc\>s_n}\>\gm_n^{-1}$.
\vvn.16>
Moreover, for any \,$m<n-1$ \,and any \,$s_1\lc s_m$ \,we have that
\,$\phi_{s_1\lc\>s_m,\<\>0\lc\<\>0}\in\pi_{m+1}(\Sg_{m+1})$ \,and
\,$\pi_{m+1}^{-1}(\phi_{s_1\lc\>s_m,\<\>0\lc\<\>0})$ \,does not depend
on \,$n$.

\vsk.2>
One can check that \,$\phi_{r_1\lc\<\>r_{n-2},\<\>0}=0$ \,if \,$r_i=0$ \,for
some \,$i=2\lc n-3$. This can be done by determining \,$\phi_{s_1\lc\>s_n}$
recursively from
\vvn.3>
\be
\exp\>\biggl(\,\sum_{s_1\lc s_n=0\!}^\infty
\phi_{s_1\lc\>s_n}\,u_1^{s_1}\?\dots\>u_n^{s_n}\biggr)\,=\,
1+\sum_{k=1}^{n-2}\,\sum_{1\le\>i_1<\<\dots<\>i_k\le n\!\!}
\vec\si_{i_1\lc i_k}\,u_{i_1}\?\dots u_{i_k}\,,
\vv.2>
\ee
see~\Ref{logn}, \,and employing induction on \,$\sum_{i=1}^{n-1} r_i$\>.
Define the elements \,$\tht_k\?\in\Sg_{k+1}$ by
\vvn.1>
\be
\pi_{k+1}(\tht_k)\,=\,
\sum_{m=1}^k\sum_{\fratop{s_1\lc\>s_m=1\!}{s_1+\dots+s_m=k\!}}^{k-m+1}
\phi_{s_1\lc\>s_m,\<\>0\lc\<\>0}\,.
\vv-.1>
\ee
Since
\vvn-.5>
\be
\log\>\Bigl(\>1+\sum_{k=1}^{n-2}\,\gm_n^{-1}\>G_{n-k}\^n\>u^k\>\Bigr)\,=\,
\sum_{k=1}^{n-2}\,I_k\^n\>u^k+\>O(u^{n-1})
\vv.4>
\ee
taking \,$u_1=\dots=u_n=u$ \,in~\Ref{logn} \,yields~\Ref{Ink}.
\end{proof}

Because of formula~\Ref{Ink}, we call the elements \;$I_1\^n\<\lc I_{n-2}\^n$
\,the {\it local charges}.

\begin{cor}
\label{generateI}
The subalgebra \,$\Aes_n$ is generated by \,$\gm_n$ and the local charges
\,$I_1\^n\<\lc I_{n-2}\^n$.
\end{cor}
\begin{proof}
The claim follows from Theorem~\ref{generateG}.
\end{proof}

Recall that \,$\Sg_n$ acts on \,$(\C^N)^{\ox n}$ by permuting the tensor
factors, and \;$\varpi_n:\C[\Sg_n]\to\End\bigl((\C^N)^{\ox n}\bigr)$ \,is
the corresponding homomorphism. Then
\vvn.2>
\be
\varpi_n(I_1\^n)\,=\,
\sum_{a=1}^n\,\sum_{i,j=1}^N\,E_{\ij}\@a\ox E_{\ji}^{(a+1)}\,,
\vv.2>
\ee
where \;$E_{\ij}\@a=\>1^{\ox(a-1)}\ox E_{\ij}\ox 1^{\ox(n-a)}$
\,and \,$E_{\ij}^{(n+1)}$ \>stands for \,$E_{\ij}\@1$.
In particular, for $N=2$ the operator
\;$n-2\>\varpi_n(I_1\^n)\in\End\bigl((\C^2)^{\ox n}\bigr)$
\,is the Hamiltonian of the {\sl XXX\/} Heisenberg model.

\appendix

\section*{Appendix}
\refstepcounter{section}

\subsection{Bethe subalgebra of \,$U(\glnt)$ and
\,Bethe algebras $\Bc_\Llb$\,}
\label{A1}
Let \,$e_{\ij}$ be the standard generators of \,$\gln$ satisfying the relations
\,$[e_{\ij}\>,e_{\kl}]=\dl_{\jk}\>e_{\il}-\dl_{\il}\>e_{\kj}$.
Let $\glnt=\gln\otimes\C[t]$ be the Lie algebra of \,$\gln$-valued polynomials
with the pointwise commutator. We identify the Lie algebra \,$\gln$ with
the subalgebra $\gln\otimes1$ of constant polynomials in $\glnt$.

\vsk.2>
Consider first-order formal differential operators in $u$\>:
\vvn.2>
\be
\Xt_{\ij}\,=\,\dl_{\ij}\,\der_u\>-
\sum_{r=0}^\infty \,(e_{\ij}\ox t^r)\,u^{-r-1}\>,\qquad i,j=1\lc N\>,
\vv.2>
\ee
and the $N$-th order formal differential operator in $u$
\vvn.1>
\be
\Dt\,=\,\sum_{\si\in\Sg_N}\sign(\si)\,
\Xt_{\si(1),1}\>\Xt_{\si(2),2}\dots\Xt_{\si(N),N}\,=\,
\der_u^N+\sum_{i=1}^N\,\sum_{j=i}^\infty\,
(-1)^i\,B_{\ij}\,u^{-j}\>\der_u^{\>N\?-i}\,.
\vv.2>
\ee
The unital subalgebra \,$\Bc\subset\Uglnt$ \,generated by the elements
\,$B_{\ij}$, \,$i=1,\dots,N$, \;$j\in\Z_{\ge i}\>$, \>is called
the {\em Bethe subalgebra\/} of \,$\Uglnt$.

\begin{thm}[\cite{T}\>, \cite{MTV1}]
The subalgebra $\Bc$ is commutative and commutes with the subalgebra
$\Ugln\subset\Uglnt$.
\qed
\end{thm}

Let \,$M$ \>be a \,$\gln$-module. Recall that \,$v\in M$ has weight
$\bs\la=(\la_1\lc\la_N)$ \,if \,$e_{\ii}\>v=\la_i\>v$ \,for all \,$i=1\lc N$,
and \,$v\in M$ is singular if \,$e_{\ij}\<\>v=0$ \,for all \,$i<j$.
\vsk.2>
Given \>$b\in\C$\>, each \,$\gln$-module \,$M$ \>becomes the evaluation module
\,$M(b)$ \>over \,$\glnt$ via the homomorphism $\glnt\to\gln$,
\;$g\ox t^r\mapsto g\>b^r$ \,for any \,$g\in\gln$, \,$r\in\Z_{ge0}$.
\vsk.2>
Let $L_{\bs\la}$ be the irreducible highest weight $\gln$-module of highest
weight $\bs\la$. For a collection $\bs\La=(\bs\la\@1\<\lc\bs\la\@k)$ of
\,$\gln$-weights, consider the tensor product \,$\ox_{i=1}^k L_{\bs\la\@i}$
of \,$\gln$-modules and denote by $\Mc_\Ll\subset\ox_{i=1}^k L_{\bs\la\@i}$
the subspace of singular vectors of weight \,$\bs\la$.
\vsk.2>
Given complex numbers \,$\bs b=(b_1\lc b_k)$, consider the tensor product
$\ox_{i=1}^k L_{\bs\la\@i}(b_i)$ of evaluation $\glnt$-modules. The action of
\>$\Bc$ on $\ox_{i=1}^k L_{\bs\la\@i}(b_i)$ preserves the subspace \,$\Mc_\Ll$.
Let \,$\Mc_\Llb$ denote the corresponding \,$\Bc$-module.
By definition, the algebras \,$\Bc_\Ll$ and \,$\Bc_\Llb$ are the images
of \,$\Bc$ in \,$\End\bigl(\ox_{i=1}^k L_{\bs\la\@i}(b_i)\bigr)$ \,and
\,$\End(\Mc_\Llb)$, respectively.
\vsk.2>
The algebras $\Bc_\Llb$ were studied in~\cite{MTV4}.
When $\bs\la\@1\<\lc\bs\la\@k$ and $\bs\la$ are partitions
with at most \>$N$ parts, and $b_1\lc b_k$ are distinct,
the $\glnt$-module $\ox_{i=1}^k L_{\bs\la\@i}(b_i)$ is irreducible.
If some of $b_1\lc b_k$ coincide, then the $\glnt$-module
$\ox_{i=1}^k L_{\bs\la\@i}(b_i)$ is a direct sum of irreducible submodules,
and these submodules are tensor products of evaluation \,$\glnt$-modules.
In such a case, the algebra $\Bc_\Llb$ is isomorphic to the direct sum
\,$\bigoplus_{\bs\La'\!,\bs b'}\Bc_{\bs\La'\!,\>\bs\la\>,\bs b'}$,
where the pairs $\bs\La'\?,\bs b'$ label nonequivalent irreducible
$\glnt$-submodules of $\ox_{i=1}^k L_{\bs\la\@i}(b_i)$.

\vsk.2>
Recall that \,$E_{\ij}\<\in\End(\C^N)$ \,is the matrix with only one nonzero
entry equal to $1$ at the intersection of the $i$-th row and $j$-th column.
The assignment \,$e_{\ij}\mapsto E_{\ij}$ \>makes \,$\C^N$ into
the \,$\gln$-module isomorphic to \,$L_{\bs\om}$\>, \,where
\,$\bs\om=(1,0\lc 0)$ is the first fundamental weight of \,$\gln$.
The algebras \,$\Bc_{n,N}\pz$ and \,$\Bc_\nnla\pz$, introduced in
Section~\ref{Gaudin}, coincide respectively with the algebras \,$\Bc_\Ll$
\>and \,$\Bc_\Llb$ \,for \,$\bs\La=(\bs\om\lc\bs\om)$ \,and \;$\bs b=\pz$.
Notice that for distinct \>$\zn$, this algebra \,$\Bc_\Llb$ is isomorphic
to the algebra \,$\Bc_{\bs\La,\bs b}$ in \cite{MTV4}, which is used in
the Theorem \ref{manyprop}.

\subsection{Algebra \,$\Oc_\Llb$}
\label{A2}
Let \,$K_1\lc K_N$ and \,$b_1\lc b_k$ \,be two collections of complex numbers,
distinct within each collection. Let \,$\bs\la\>,\>\bs\la\@1\<\lc\>\bs\la\@k$
\,be partitions with at most \>$N$ parts. The algebra \,$\Oc_\Llb$, \,where
\,$\bs\La=(\bs\la\@1\<\lc\bs\la\@k)$\>, \;$\bs b=(b_1\lc b_k)$ \,and
the dependence on \,$K_1\lc K_N$ \,is suppressed, was defined for these data
and studied in~\cite[Section~5]{MTV2}.
\vsk.2>
The algebra \,$\Ee_{\bs\mu}\pz$, used in the proof of
Theorem~\ref{BSnH}, coincides with the algebra $\Oc_\Llb$ \,for the following
choice of parameters: \,$N=n$\>, \;$k=1$\>, \;$b_1=0$\>,
\;$\bs\la\@1\<=\bs\mu$\>, \;$\bs\la=(1\lc 1)$\>, \,and \,$K_i=z_i$
\,for all \,$i=1\lc n$.

\subsection{Proof of Lemma~\ref{PePe}}
\label{A4}
In the proof, we are using several identities for rational functions whose
verification is left to a reader.

\vsk.3>
Let \,$C$ \>be the matrix with entries
\vvn.2>
\;$C_{ab}=\prod_{c=1}^{a-1}\,(z_b-z_c)$ \;for \,$a\le b$\>, \,and
\,$C_{ab}=0$ \,for \,$a>b$\>. Then the entries of \,$C^{-1}$ are
\,$(C^{-1})_{ab}=\prod_{c=1,\,c\ne a}^b\,(z_a-z_c)^{-1}$
for \,$a\le b$\>, \,and \,$(C^{-1})_{ab}=0$ \,for \,$a>b$\>.

\vsk.2>
Recall that \,$Z$ \,and \,$\Zh$ \,are the matrices with entries
\,$Z_{ab}=z_a\>\dl_{ab}$ \,and \,$\Zh_{ab}=\dl_{a,b-1}$.
It is straightforward to see that
\beq
\label{CZC}
CZ\>C^{-1}\,=\,Z+\Zh\,.
\eeq
\par
Let \;$Q_\hbar$ \>be the matrix given by \Ref{Qh}, \Ref{caz}, \,and
\;$\Qt=CQ_\hbar C^{-1}$ \,with \,$\hbar=1$. Then after a minor simplification
\be
\Qt_{ab}\>=\,\sum_{c=a}^n\,\sum_{d=1}^b\,\frac{q(z_c)}{z_c-z_d+1}\;
\prod_{\fratop{r=a}{r\ne c}}^n\,\frac1{z_c-z_r}\;
\prod_{\fratop{s=1}{s\ne d}}^b\,\frac1{z_d-z_s}\;.
\vv-.4>
\ee
Taking the sum over \,$d$\>, we get
\vvn.1>
\be
\Qt_{ab}\>=\,\sum_{c=a}^n\,q(z_c)\;
\prod_{\fratop{r=a}{r\ne c}}^n\,\frac1{z_c-z_r}\,
\prod_{s=1}^b\,\frac1{z_c-z_s+1}\;=\>
\frac1{2\<\>\pi\<\>i}\int\limits_{\!\!\gm}q(u)\;
\prod_{r=a}^n\,\frac1{u-z_r}\,\prod_{s=1}^b\,\frac1{u-z_s+1}\;du\,,
\vvn.2>
\ee
where \,$\gm$ \>is a simple closed curve such that the points \,$z_1\lc z_n$
\>are inside \,$\gm$ \>and \,$z_1-1\lc z_n-1$ \>are outside. For instance,
if \,$|\>z_a|<1/3$ \,for all \,$a$\>, one can take \,$\gm$ \>to be the circle
\,$|\>u\>|=1/2$\>. Therefore, in the limit \;$z_a\?\to 0$ \;for all
\,$a=1\lc n$\>,
\vvn.4>
\be
\Qt_{ab}\,\to\,\frac1{2\<\>\pi\<\>i}\int\limits_{\!\!|u|=1\</2}
\!\frac{q(u)}{u^{n-\<\>a\<\>+1}\>(u+1)^{\<\>b}}\;du\,=\,\Qh_{ab}\,,
\vv.3>
\ee
cf.~\Ref{Qh}. The last formula together with \Ref{CZC}, \Ref{Peh} and
\Ref{Peht} implies that
\vvn.4>
\be
\Pe_\hbar(u,v;\zn;q)\,\to\,\Peh(u,v;q)\,,\qquad z_a\?\to0\,,\quad a=1\lc n\,,
\vv.3>
\ee
which proves Lemma~\Ref{PePe}.
\qed

\end{document}